\documentclass{amsart}

\usepackage{color}
\usepackage[dvips]{graphicx}
\usepackage{amscd}
\usepackage{amsmath}
\usepackage{ascmac}
\usepackage{cases}
\usepackage{amssymb}
\usepackage{amsfonts}
\usepackage{amsthm}
\usepackage[all]{xy}
\usepackage{enumerate}
\usepackage{wrapfig}

\theoremstyle{definition}
\newtheorem{Thm}{{\bf Theorem}}[section]

\newtheorem{Prop}[Thm]{{\bf Proposition}}

\newtheorem{Def}[Thm]{{\bf Definition}}
\newtheorem{Rem}[Thm]{{\bf Remark}}

\newcounter{Exami}

\newcommand{\bnu}{{\boldsymbol \nu}}

\newcommand{\bxi}{{\boldsymbol \xi}}
\newcommand{\bt}{{\bold t}}

\newcommand{\bC}{{\mathbb C}}       
\newcommand{\bP}{{\mathbb P}}      
\newcommand{\bZ}{{\mathbb Z}}

\newcommand{\cO}{{\mathcal O}}

\newcommand{\cM}{{\mathcal M}}

\newcommand{\cI}{{\mathcal I}}

\newcommand{\cK}{{\mathcal K}}

\newcommand{\cV}{{\mathcal V}}

\newcommand{\Aut}{\mathop{\rm Aut}\nolimits}

\newcommand\Hom{\mathop{\rm Hom}\nolimits}

\newcommand\Spec{\mathop{\rm Spec}\nolimits}
\newcommand\Hilb{\mathop{\rm Hilb}\nolimits}

\newcommand{\res}{\mathop{\sf res}\nolimits}

\newcommand\lra{\longrightarrow}
\newcommand\ra{\rightarrow}

\title[Explicit description of jumping phenomena]%
  {Explicit description of jumping phenomena on moduli spaces of parabolic connections and Hilbert schemes of points on surfaces}
\author[A.\ Komyo and M.-H.\ Saito]{Arata Komyo \and Masa-Hiko Saito}
\address{Department of Mathematics, Graduate School of Science, Kobe University, 1-1 Rokkodai-cho, Nada-ku, Kobe, 657-8501, Japan}
\email{akomyo@math.kobe-u.ac.jp}
\address{Department of Mathematics, Graduate School of Science, Kobe University, 1-1 Rokkodai-cho, Nada-ku, Kobe, 657-8501, Japan}
\email{mhsaito@math.kobe-u.ac.jp}
\subjclass[2010]{Primary~14D20, Secondary~34M55~32G34.}
\keywords{Parabolic connection, Parabolic Higgs bundle, Apparent singularity.}
\thanks{ This research was partly supported by JSPS Grant-in-Aid for  
Scientific Research (S)24224001, challenging Exploratory Research 15K13427}

\begin{document}

\maketitle

\begin{abstract}
In this paper, we investigate the apparent singularities and the dual parameters of 
rank $2$ parabolic connections on $\bP^1$ and
rank $2$ (parabolic) Higgs bundle on $\bP^1$.
Then we obtain explicit descriptions of Zariski open sets of the moduli space of the parabolic connections and the moduli space of the Higgs bundles.
For $n=5$, we can give global descriptions of the moduli spaces in detail.
\end{abstract}

\section{Introduction}

The purpose of this paper is to give explicit descriptions of the moduli space of rank $2$ (parabolic) Higgs bundles on $\mathbb{P}^1$ and 
the moduli space of rank $2$ parabolic connections on $\mathbb{P}^1$ by the \textit{apparent singularities} and their \textit{dual parameters}.
It is well-known that the apparent singularities and their dual parameters are coordinates on Zariski open sets of these moduli spaces.
Historically,
Okamoto \cite{Okamoto} described Hamiltonians of the Garnier systems by the apparent singularities and their duals,
for the Garnier systems are obtained by the isomonodromic deformations of rank $2$ parabolic connections on $\mathbb{P}^1$.
The apparent singularities and their duals are introduced as coordinates for a Zariski open set of Okamoto's space of initial conditions,
which are nothing but the moduli space of rank $2$ parabolic connections on $\mathbb{P}^1$.
Arinkin--Lysenko \cite{AL} and Oblezin \cite{Obl} studied the moduli space of rank $2$ parabolic connections on $\mathbb{P}^1$ more systematically,
and they also introduced the apparent singularities and their duals as coordinates for Zariski open sets of moduli spaces.
Oblezin also showed that the moduli space of rank $2$ parabolic connections on $\mathbb{P}^1$ is birational to the Hilbert scheme of points on the blowing up of the total space
of a certain line bundle on $\mathbb{P}^1$.
Dubrovin and Mazzoco discussed the apparent singularities and their duals for higher rank cases on $\mathbb{P}^1$ in detail \cite{DM}.
For the moduli spaces of parabolic connections on arbitrary genus curves, Inaba--Iwasaki--Saito \cite{IIS} and Inaba \cite{Inaba} established the existence of good moduli spaces
of stable parabolic connections,
and it is interesting to describe their geometric structures.
Saito and Szabo are developing a systematic treatment of apparent singularities and their duals for general parabolic connections on higher genus curves \cite{SS}.
One can show the similar geometric description of the moduli spaces of parabolic Higgs bundles.
The main purpose of this paper is to give more explicit description of the total space of the moduli spaces of rank $2$ parabolic connections on $\mathbb{P}^1$.
For the purpose, we need treat the following particular cases.
The first case is that the apparent singularities approaches to the regular singularities of Higgs fields or connections 
(This case is already treated in \cite[Section 3.7]{Obl}).
The second case is that the apparent singularities have multiplicities.
The third case is that the type of the underlying bundle is jumping.
The jumping phenomenon happens to $n$-regular singularities cases where $n\ge 5$. 
For the first and second cases, we can give an explicit description of families for the $n$-point regular singularities case.
For the third case, we give an explicit description of jumping families parameterized by the apparent singularities and their duals for the $5$-point regular singularities case.
Oblezin \cite{Obl} considered the stratification of the moduli spaces of rank $2$ parabolic connections on $\mathbb{P}^1$ associated to the bundle type of underlying bundles,
and gave geometric descriptions of each strata of the moduli spaces, separately.
On the other hand, in this paper, we try to give a global geometric description of the moduli spaces including the jumping phenomena of bundle type.
As the result, we obtain a global description of the moduli spaces for the $5$-point regular singularities case,
and give an explicit description of universal families of Higgs bundles and connections.

Fix points $t_1,\ldots, t_n \in \bP^1$ $(t_i\neq t_j)$, and set $D=t_1+\cdots +t_n$.
We consider pairs $(E,\nabla)$ where $E$ is a rank $2$ vector bundle on $\bP^1$ and $\nabla\colon E \ra E\otimes \Omega^1_{\bP^1}(D)$
a connection having simple poles supported on $D$.
At each pole, we have two residual eigenvalues $\{ \xi^+_i, \xi^-_i\}$ of $\nabla$, $i=1,\ldots,n$;
they satisfy Fuchs relation $\sum_i(\xi^+_i+\xi^-_i)+d=0$ where $d=\deg (E)$.
Moreover, we introduce parabolic structures $\boldsymbol{l}=\{l_i\}_{1\le i \le n}$
such that $l_i$ is a one dimensional subspace of $E|_{t_i}$ which corresponds to an eigenspace of the residue of $\nabla$ at $t_i$ with the eigenvalue $\xi^+_i$.
Note that when $\xi^+_i \neq \xi^-_i$, the parabolic structure $\boldsymbol{l}$ is determined uniquely by the connection $(E,\nabla)$.
Fixing a spectral data $\bxi=(\xi^{\pm}_i)$ with integral sum $-d$ and introducing the weight $\boldsymbol{w}$,
we can construct the moduli space $M^{\boldsymbol{w}}_{\bt, \bxi}(\mathfrak{gl}_2)$ of \textit{$\boldsymbol{w}$-stable $\bxi$-parabolic connections} $(E,\nabla, \boldsymbol{l})$
by Geometric Invariant Theory and 
the moduli space $M^{\boldsymbol{w}}_{\bt, \bxi}(\mathfrak{gl}_2)$ turns to be a smooth irreducible quasi-projective variety of dimension $2n-6$
for generic weight $\boldsymbol{w}$ (see \cite{IIS}).
Note that, when
\begin{equation}\label{generic condition 1}
\sum^n_{i=1}  \xi_i^{\epsilon_i} \notin \bZ
\end{equation}
for any $(\epsilon_i), \epsilon_i \in \{ +,- \}$,
every parabolic connection $(E,\nabla,\boldsymbol{l})$ is irreducible, hence stable.
Therefore the moduli space $M^{\boldsymbol{w}}_{\bt, \bxi}(\mathfrak{gl}_2)$ does not depend on the choice of the weight $\boldsymbol{w}$ in such cases.
It is known that the moduli spaces coincide with the spaces of initial conditions for Garnier systems,
and the case $n=4$ corresponding to the Pinlev\'e VI equation, 
for such differential equations are nothing but isomonodromic deformations for linear connections.
Next, we fix $\bxi=(\xi^{\pm}_i)_{1\le i \le n}$ where $\sum_i(\xi^+_i+\xi^-_i)=0$. 
In the same way as above, 
we can define \textit{$\boldsymbol{w}$-stable $\bxi$-parabolic Higgs bundle $(E,\Phi, \boldsymbol{l})$}.
Here $E$ is a rank $2$ vector bundle on $\bP^1$, $\Phi\colon E \ra E\otimes \Omega^1_{\bP^1}(D)$ is an $\cO_{\mathbb{P}^1}$-morphism, 
and $\boldsymbol{l}$ is the parabolic structure.
At each point $t_i$, residual eigenvalues of $\Phi$ are $\{ \xi^+_i, \xi^-_i\}$ of $\Phi$.
Let $M^{\boldsymbol{w}}_{H, \bt, \bxi}(\mathfrak{gl}_2)$ be the moduli space of $\boldsymbol{w}$-stable $\bxi$-parabolic Higgs bundles.

By suitable transformations, we may assume that $d=\deg (E)= -1$ and $\bxi$ can be normalized as follows.
For connection cases, we can put
\begin{equation*}
\begin{cases}
\xi_i^+=\nu_i& (i=1,\ldots,n)\\
\xi_i^-=- \nu_i& (i=1,\ldots,n-1) \\
\xi_n^-=1 - \nu_n,
\end{cases}
\end{equation*}
and for Higgs cases, we can put $\xi_i^+=\nu_i$ and $\xi_i^-=- \nu_i$ ($i=1,\ldots,n$),
for some $\bnu=(\nu_1,\ldots,\nu_n) \in \mathbb{C}^{n}$.
Let $M$ and $M_H$ be the moduli space of $\bnu$-$\mathfrak{sl}_2$-parabolic connections 
and the moduli space of $\bnu$-$\mathfrak{sl}_2$-parabolic Higgs bundles, respectively.
By these normalizations, we have natural isomorphisms $M \cong M^{\boldsymbol{w}}_{\bt, \bxi}(\mathfrak{gl}_2)$ and 
$M_H \cong M^{\boldsymbol{w}}_{H, \bt, \bxi}(\mathfrak{gl}_2)$.
(Note that the moduli space $M$ is noting but the moduli space of (modified) $(\nu_1,\ldots,\nu_n)$-bundles on $\mathbb{P}^1$ treated in \cite{AL} and \cite{Obl}.)

For the moduli space $M_{{\it H}}$, we obtain the following results.
First, we consider the Zariski open set $M_{{\it H}}^0$ which is the locus where the type of the underlying bundle is $\cO\oplus \cO(-1)$.
Let $\cK'_n$ be some Zariski open set of some blowing-up of the Hirzebruch surface of degree $n-2$. (See Figure \ref{fig:one}).
By the explicit computation of the apparent singularities and their dual parameters, we have the following 
\begin{Thm}[Theorem \ref{Prop M^0 to Hilb}]\label{Prop M^0 to Hilb 0}
\textit{
By the apparent singularities and dual parameters, we have a map 
\begin{equation*}
M_{{\it H}}^0 \lra \mathrm{Hilb}^{n-3}(\cK'_n),
\end{equation*}
and this map is injective.
Moreover, we can give an explicit description of the universal family $(\tilde{E}^{(0)},\tilde{\Phi}^{(0)}) \rightarrow M_{{\it H}}^0\times \mathbb{P}^1$.
}
\end{Thm}

Suppose $n=5$. 
We consider the total moduli space $M_{{\it H}}$, which also includes the \textit{jumping locus}. 
The type of the underlying bundle of members are $\cO\oplus \cO(-1)$ (generic) or $\cO(1)\oplus \cO(-2)$ (jumping locus).
Let $\widehat{M}_{{\it H}}$ be the moduli space of \textit{$\bnu$-$\mathfrak{sl}_2$-parabolic Higgs bundles with a cyclic vector $\sigma \in H^0(\mathbb{P}^1, E)$}.
The moduli space $\widehat{M}_{{\it H}}$ is the blowing-up of $M_{\it H}$ along the jumping locus.
We take some blowing-up of $\mathrm{Hilb}^2(\cK'_5)$, denoted by $\widetilde{\mathrm{Hilb}}^2(\cK'_5)$.
Then we have the map $\widehat{M}_{{\it H}} \ra \widetilde{\mathrm{Hilb}}^2(\cK'_5)$.
\begin{Thm}[Theorem \ref{main theorem +} and Section \ref{construction of the family on cV1}]\label{main theorem}
\textit{
Suppose that $n=5$.
The map
\begin{equation*}
\widehat{M}_{{\it H}} \lra \widetilde{\mathrm{Hilb}}^2(\cK'_5)
\end{equation*}
is injective.
Moreover, we can give an explicit description of the universal family $(\tilde{E},\tilde{\Phi}, \tilde{\sigma}) \rightarrow \widehat{M}_{{\it H}}\times \mathbb{P}^1$.
}
\end{Thm}

For the moduli space $M$ of connections, which is isomorphic to $M^{\boldsymbol{w}}_{\bt, \bxi}(\mathfrak{gl}_2)$, we have the following results.
First, we consider the Zariski open set $M^0$ which is the locus where the type of the underlying bundle is $\cO\oplus \cO(-1)$.
Let $\widetilde{\cK}_n'$ be some Zariski open set of some blowing-up of the Hirzebruch surface of degree $n-2$.
By the same argument as in the Higgs case, we have the following 
\begin{Thm}[Theorem \ref{from M0con to Hilb}]\label{from M0con to Hilb 0}
\textit{
By the apparent singularities and dual parameters, we have a map 
\begin{equation}\label{M0 to Hilb 0}
M^0 \lra \mathrm{Hilb}^{n-3}(\widetilde{\mathcal{K}}_n'),
\end{equation}
and this map is injective.
Moreover, we can give an explicit description of the universal family $(\tilde{E}^{(0)},\tilde{\nabla}^{(0)}) \rightarrow M^0\times \mathbb{P}^1$.
}
\end{Thm}

Parts of Theorem \ref{Prop M^0 to Hilb 0} and Theorem \ref{from M0con to Hilb 0} are already contained in \cite{AL}, \cite{IIS2}, \cite{Obl}, and \cite{Yoshi}.
For $n=4$, the results are discussed in \cite{AL} and \cite{IIS2}. 
Oblezin \cite{Obl} gives a map from $M_{{\it H}}^0$ (resp. $M^0$) to a $(n-3)$-th symmetric product of $\cK'_n$ (resp. $\widetilde{\mathcal{K}}_n'$),
which is an isomorphism on a certain open set.
For $n=5$, the injectivities are discussed in \cite{Yoshi}.

Suppose $n=5$. 
We consider the moduli space $M$, which includes the jumping locus. 
The type of the underlying bundle of members of the jumping locus is $\cO(1)\oplus \cO(-2)$.
Let $\widehat{M}$ be the moduli space of \textit{$\bnu$-$\mathfrak{sl}_2$-parabolic connections with a cyclic vector $\sigma \in H^0(\mathbb{P}^1, E)$}.
The moduli space $\widehat{M}$ is the blowing-up of $M$ along the jumping locus.
Let $C_{\infty}$ be the $\infty$-section of the Hirzebruch surface of degree $n-2$.

\begin{Thm}[Theorem \ref{main theorem conn}]\label{main theorem conn intro}
\textit{
Let $\phi \colon \widehat{M} \dashrightarrow \Hilb^2(\widetilde{\cK}_n'\cup C_{\infty})$ be the birational map constructed by the apparent singularities and the dual parameters.
By taking some sequence of blowing-ups $\widetilde{\Hilb}^2(\widetilde{\cK}_n'\cup C_{\infty}) \ra \Hilb^2(\widetilde{\cK}_n'\cup C_{\infty})$,
we have the injective map $\tilde{\phi} \colon \widehat{M} \rightarrow \widetilde{\Hilb}^2(\widetilde{\cK}_n'\cup C_{\infty})$ for $\phi$.
The moduli space $\widehat{M}$ is biregular to its image $\mathrm{Image}\, \tilde{\phi}(\widehat{M}) \subset \widetilde{\Hilb}^2(\widetilde{\cK}_n'\cup C_{\infty})$.
}
\end{Thm}

The organization of this paper is as follows.
In Section 2, we introduce definitions and notations which are necessary in this paper.
In Section 3, \ref{O+O(-1)}, we consider the moduli spaces of $\bnu$-$\mathfrak{sl}_2$-parabolic Higgs bundles with bundle type $\cO \oplus \cO(-1)$.
We show Theorem \ref{Prop M^0 to Hilb 0} (Theorem \ref{Prop M^0 to Hilb}) by explicit calculations of apparent singularities.
In \ref{O(1)+O(-2)}, we consider the moduli spaces of $\bnu$-$\mathfrak{sl}_2$-parabolic Higgs bundles with a cyclic vector for $n=5$.
We show the injectivity of the map in Theorem \ref{main theorem} (Theorem \ref{main theorem +}) 
by explicit calculations of apparent singularities and spectral curves.
In Section 4, we construct an explicit jumping families of $\bnu$-$\mathfrak{sl}_2$-parabolic Higgs bundles by the lower and upper modifications.
In particular, in \ref{construction of the family on cV1}, we give an explicit description of the universal family of $\widehat{M}_{{\it H}}$
In Section 5, we consider the moduli spaces of $\bnu$-$\mathfrak{sl}_2$-parabolic connections.
In \ref{connection O+O(-1)}, we show Theorem \ref{from M0con to Hilb 0} (Theorem \ref{from M0con to Hilb}) by the same way as in the Higgs case.
In \ref{jumping phenomenon conn}, we construct an explicit jumping family of connections for $n=5$, and
in \ref{lim of app and dual}, we analyze the behavior of the apparent singularities and their duals when the parameters of the jumping family approach to the jumping locus.
Finally, we obtain Theorem \ref{main theorem conn intro} (Theorem \ref{main theorem conn}).

\section{Preliminaries}\label{Preliminaries}

In this section, first, we define \textit{$\bnu$-$\mathfrak{sl}_2$-parabolic connections} and \textit{$\bnu$-$\mathfrak{sl}_2$-parabolic Higgs bundles}, 
and recall the well-known facts of the connections and Higgs bundles.
In 2.2, we describe some blowing-ups of the Hirzebruch surface $\Sigma_{n-2}$, which are target spaces of the map defined by the apparent singularities and their duals.
In 2.3, we discuss descriptions of Higgs fields.
Since we consider Higgs bundles on $\mathbb{P}^1$, the underlying vector bundles split into the direct sum of line bundles.
Then we can describe the Higgs fields explicitly.
By the automorphisms of vector bundles, we can normalize the Higgs fields to reduce the number of parameters.
In 2.4, we discuss the apparent singularities and the dual parameters of Higgs bundles, which
give a map from the moduli space of Higgs bundles (with a cyclic vector) to the symmetric product of the total space of $\Omega^1_{\mathbb{P}^1}(D)$.
In 2.5, we discuss the transformations called \textit{the lower modification} and \textit{the upper modification},
and we use the transformations for a construction of a universal family of the moduli space of Higgs bundles (with a cyclic vector).
The contents of 2.1, 2.2, and 2.5 basically follow the expositions from \cite{AL} and \cite{Obl}.

\subsection{$\mathfrak{sl}_2$-connections and $\mathfrak{sl}_2$-Higgs bundles}

We introduce $\mathfrak{sl}_2$-parabolic connections and $\mathfrak{sl}_2$-parabolic Higgs bundles,
and we consider relations between the moduli space $M^{\boldsymbol{w}}_{\bt, \bxi}(\mathfrak{gl}_2)$ and these moduli spaces.

Fix complex numbers $\nu_1,\ldots,\nu_n \in \bC$.
Suppose that $\nu_1 \cdots \nu_n \neq 0$ and 
\begin{equation}\label{generic condition 2}
\sum^n_{i=1} \epsilon_i \nu_i \notin \bZ
\end{equation}
for any $(\epsilon_i), \epsilon_i \in \{ 1,-1 \}$.
\begin{Def}
A \textit{$\bnu$-$\mathfrak{sl}_2$-parabolic connection on $\bP^1$} is a triplet $(E, \nabla,\varphi)$
such that 
\begin{enumerate}
\item $E$ is a rank $2$ vector bundle on $\bP^1$,
\item $\nabla\colon E \ra E\otimes \Omega^1_{\bP^1}(D)$ is a connection,
\item $\varphi\colon \bigwedge^2 E \cong \cO_{\bP^1}(-1)$ is a horizontal isomorphism, 
\item the residue $\res_{t_i}(\nabla)$ of the connection $\nabla$ at $t_i$ has eigenvalues $\nu_i^{\pm}$, $1\le i \le n$.
\end{enumerate}
Here we put 
\begin{equation*}
 \nu^{\pm}_i := \pm \nu_i \ \ (i=1,\ldots, n-1 ),\ \nu^+_n:=\nu_n,\ \nu^-_n := 1-\nu_n.
\end{equation*}
\end{Def}
Denote by $\cM$ the moduli stack of $\bnu$-$\mathfrak{sl}_2$-parabolic connections on $\bP^1$,
and by $M$ its coarse moduli space.
Let $\nabla' \colon \cO_{\mathbb{P}^1} \ra \cO_{\mathbb{P}^1} \otimes \Omega^1_{\mathbb{P}^1}(D)$ be the connection defined by
\begin{equation*}
f \longmapsto df + \left( \frac{1}{2} \sum^{n-1}_{i=1}  (-\xi^+_i - \xi_i^-) \frac{dz}{z-t_i} + \frac{1}{2} \sum^{n-1}_{i=1}  (\xi^+_i + \xi_i^-) \frac{dz}{z-t_n} \right) f.
\end{equation*}
Suppose that the condition (\ref{generic condition 1}) holds and $\xi_i^{+}\neq \xi_n^{-}$ for $i=1,\ldots,n$. 
Then we have an isomorphism
\begin{equation*}
\begin{aligned}
 M^{\boldsymbol{w}}_{\bt, \bxi}(\mathfrak{gl}_2) &\lra M = M_{\bnu} \\
 (E,\nabla ,\boldsymbol{l}) &\longmapsto ((E, \nabla) \otimes (\cO_{\mathbb{P}^1}, \nabla'), \varphi),
\end{aligned}
\end{equation*}
where $\bnu=(\nu_1,\ldots, \nu_n)$, $\nu_i= (\xi_i^{+}-\xi_i^{-})/2$ ($i=1,\ldots, n-1$), $\nu_n=\xi^+_n+\sum_{i=1}^{n-1}(\xi_i^{+}-\xi_i^{-})/2$.
\begin{Def}
A \textit{$\bnu$-$\mathfrak{sl}_2$-parabolic Higgs bundle on $\bP^1$} is a triplet $(E, \Phi,\varphi)$
such that 
\begin{enumerate}
\item $E$ is a rank $2$ vector bundle on $\bP^1$,
\item $\Phi \colon E \ra E\otimes \Omega^1_{\bP^1}(D)$ is an $\cO_{\bP^1}$-morphism,
\item $\varphi\colon \bigwedge^2 E \cong \cO_{\bP^1}(-1)$ is an isomorphism and $\mathrm{tr}(\Phi)=0$,
\item the residue $\res_{t_i}(\Phi)$ of the connection $\Phi$ at $t_i$ has eigenvalues $\pm \nu_i$, $1\le i \le n$.
\end{enumerate}
\end{Def}
Denote by $\cM_{{\it H}}$ the moduli stack of $\bnu$-$\mathfrak{sl}_2$-parabolic Higgs bundles on $\bP^1$,
and by $M_{{\it H}}$ its coarse moduli space.
We have a stratification of $M_{{\it H}}$ as follows.
By the irreducibility of $(E ,\Phi ,\varphi) \in M_{{\it H}}$, we have the following 
\begin{Prop}\label{bundle type}
\textit{
For $(E ,\Phi ,\varphi) \in M_{{\it H}}$,
we have 
\begin{equation*}
E\cong \cO(k) \oplus\cO(-k-1) \text{  where } 0\le k \le \left[ \frac{n-3}{2} \right].
\end{equation*}
}
\end{Prop}
Let $M^k_{{\it H}}$ be the subvariety of $M_{{\it H}}$ where $E\cong \cO(k) \oplus\cO(-k-1)$.
Then
\begin{equation*}
M_{{\it H}} = M_{{\it H}}^0\cup \cdots \cup M_{{\it H}}^{[ (n-3)/2]}.
\end{equation*}
Note that the stratum $M_{{\it H}}^0$ is a Zariski open dense of $M_{{\it H}}$.

Moreover, we introduce \textit{$\bnu$-$\mathfrak{sl}_2$-parabolic connection on $\bP^1$ with a cyclic vector} 
and \textit{$\bnu$-$\mathfrak{sl}_2$-parabolic Higgs bundle on $\bP^1$ with a cyclic vector}. 
\begin{Def}
A \textit{$\bnu$-$\mathfrak{sl}_2$-parabolic connection on $\bP^1$ with a cyclic vector} is a tuple $(E, \nabla,\varphi,[\sigma])$
such that 
\begin{enumerate}
\item $E$ is a rank $2$ vector bundle on $\bP^1$,
\item $\nabla \colon E \ra E\otimes \Omega^1_{\bP^1}(D)$ is a connection,
\item $\varphi\colon \bigwedge^2 E \cong \cO_{\bP^1}(-1)$ is a horizontal isomorphism,
\item the residue $\res_{t_i}(\nabla)$ of the connection $\nabla$ at $t_i$ has eigenvalues $\nu_i^{\pm}$, $1\le i \le n$.
\item $[\sigma] \subset H^0(\bP^1, E) $ is a $1$-dimensional subspace generated by a nonzero section $\sigma \in H^0(\bP^1, E)$.
\end{enumerate}
\end{Def}
Denote by $\widehat{\cM}$ the moduli stack of $\bnu$-$\mathfrak{sl}_2$-parabolic connections on $\bP^1$ with a cyclic vector,
and by $\widehat{M}$ its coarse moduli space.
Since $\dim H^0(\mathbb{P}^1,\cO\oplus \cO(-1))=1$, $M^0$ is contained in $\widehat{M}$.
\begin{Def}
A \textit{$\bnu$-$\mathfrak{sl}_2$-parabolic Higgs bundle on $\bP^1$ with a cyclic vector} is a tuple $(E, \Phi,\varphi,[\sigma])$
such that 
\begin{enumerate}
\item $E$ is a rank $2$ vector bundle on $\bP^1$,
\item $\Phi \colon E \ra E\otimes \Omega^1_{\bP^1}(D)$ is an $\cO_{\bP^1}$-morphism,
\item $\varphi\colon \bigwedge^2 E \cong \cO_{\bP^1}(-1)$ is an isomorphism and $\mathrm{tr}(\Phi)=0$,
\item the residue $\res_{t_i}(\nabla)$ of the connection $\nabla$ at $t_i$ has eigenvalues $\pm \nu_i$, $1\le i \le n$.
\item $[\sigma] \subset H^0(\bP^1, E) $ is a $1$-dimensional subspace generated by a nonzero section $\sigma \in H^0(\bP^1, E)$.
\end{enumerate}
\end{Def}
Denote by $\widehat{\cM}_{{\it H}}$ the moduli stack of $\bnu$-$\mathfrak{sl}_2$-parabolic Higgs bundles on $\bP^1$ with a cyclic vector,
and by $\widehat{M}_{{\it H}}$ its coarse moduli space.
For $n=5$, the moduli space $\widehat{M}_{{\it H}}$ is the blowing-up of $M_{{\it H}}$ along $M^1_{{\it H}}$.

\subsection{Hirzebruch surfaces and the blowing-ups}\label{Hirzebruch blow up}

For description of the moduli spaces $M$ and $M_{{\it H}}$, we introduce some blowing-ups of the Hirzebruch surface $\Sigma_{n-2}$.
Put $L:=\Omega^1_{\bP^1}(D)$.
Let $\mathbb{L}$ be the total space of the line bundle $L$.
Note that $\mathbb{L} = \Sigma_{n-2} \setminus C_{\infty}$
where $C_{\infty}$ is the infinity section $(C_{\infty})^2=-(n-2)$.

\begin{figure}[htbp]
 \begin{center}
  \includegraphics[width=150mm]{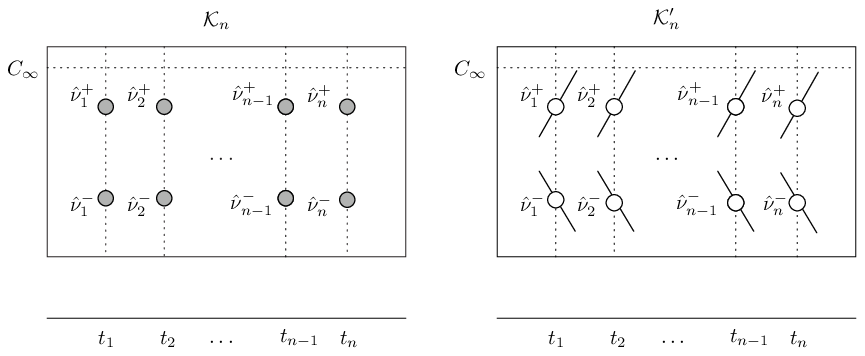}
 \end{center}
 \caption{$\cK_n$ and $\cK_n'$}
 \label{fig:one}
\end{figure}

First, we construct a blowing-up of the Hirzebruch surface $\Sigma_{n-2}$ corresponding to $M_H$.
Let $\pi \colon \mathbb{L} \ra \bP^1$ be the projection and let $\tau_i \colon \pi^{-1}(t_i) \xrightarrow{\cong} \bC$ be the residue map.
Put $\nu^{+}_i:=\nu_i$, $\nu^{-}_i:=-\nu_i$ for $i=1,\ldots,n$, and $\hat{\nu}^{\pm}_i := \tau^{-1}_i(\nu^{\pm}_i)$.
Set 
\begin{equation*}
\cK_n' := \left( \mathrm{Bl}_{\hat{\nu}^{\pm}_i} \mathbb{L} \right) \setminus ( \widetilde{F}_1\cup \cdots \cup \widetilde{F}_n ) 
\end{equation*}
where $\mathrm{Bl}_{\hat{\nu}^{\pm}_i} \mathbb{L}$ is the blowing-up of $\mathbb{L}$ at $\hat{\nu}^{\pm}_i$ for $i=1,\dots,n$,
and $\widetilde{F}_i$ are the proper pre-images of the fiber $F_i$ $i=1,\ldots,n$.
We denote by $\cK_n$ the image of $\cK_n'$ by the projection $\cK_n' \ra \mathbb{L}$ (see Figure \ref{fig:one}).

Second, we construct a blowing-up of the Hirzebruch surface $\Sigma_{n-2}$ corresponding to $M$.
Let $\pi \colon \mathbb{L}  \ra \bP^1$ be the projection and let $\tau_i \colon \pi^{-1}(t_i) \xrightarrow{\cong} \bC$ be the residue map.
Set
\begin{equation*}
\widetilde{\cK}_n' := \left( \mathrm{Bl}_{\hat{\nu}^{\pm}_i} \mathbb{L} \right) \setminus ( \widetilde{F}_1\cup \cdots \cup \widetilde{F}_n ) 
\end{equation*}
where $\hat{\nu}^{\pm}_i := \tau^{-1}_i(\nu^{\pm}_i)$. 
Here, $\nu^{+}_i:=\nu_i$, $\nu^{-}_i:=-\nu_i$ for $i=1,\ldots,n-1$ and $\nu^{+}_n:=\nu_n$, $\nu^{-}_n:=1-\nu_n$.
We denote by $\widetilde{\cK}_n$ the image of $\widetilde{\cK}_n'$ by the projection $\widetilde{\cK}_n' \ra \mathbb{L}$.

\subsection{Description of Higgs fields}

Put $U_0 := \bP^1\setminus \{\infty \}$, $U_{\infty} := \bP^1\setminus \{0 \}$.
Let $z$ and $w$ be the coordinates on $U_0$ and $U_{\infty}$, respectively.
Put 
\begin{equation}\label{omega and R}
\omega_z:=\frac{dz}{z(z-1)(z-x_1)\cdots(z-x_{n-3})}\ \text{ and }\
R_k:=\left(
\begin{array}{clcl}
z^k & 0 \\
0 & \frac{1}{z^{k+1}}
\end{array}
\right) ,\ 0\le k \le \left[ \frac{n-3}{2} \right].
\end{equation}
We consider an explicit description of the Higgs field of $(E,\Phi,\varphi) \in M_{{\it H}}$.
Suppose that $E \cong \cO(k)\oplus \cO(-k-1)$ where $0\le k \le \left[ (n-3)/2 \right]$.
We can describe the Higgs field $\Phi$ as follows:
\begin{equation}\label{normal form of O+O(-1)}
\Phi = 
\begin{cases}
 A_z^k \otimes \omega_z    & \text{on } U_0 \\
 R_k^{-1}(A_z^k \otimes \omega_z )R_k  & \text{on } U_{\infty},\quad 
\end{cases}
A_z^k :=\left(
\begin{array}{clcl}
f_{11}^{(n-2)}(z) & f_{12}^{(n+2k-1)}(z) \\
 f_{21}^{(n-2k-3)}(z) & -f_{11}^{(n-2)}(z)
\end{array}
\right)
\end{equation}
where $f_{ij}^{(l)}(z)$ is a polynomial in $z$ of degree at most $l$.
By the irreducibility, we have $f_{21}^{(n-2k-3)}(z)\neq0$.

We consider automorphisms of the vector bundle $E\cong \cO(k)\oplus \cO(-k-1)$.
Any element 
$P\in \Hom_{\cO_{\bP^1}} (E,E) \cong H^0(\bP^1, \cO\oplus \cO(2k+1) \oplus \cO)$ 
is described as follows:
\begin{equation*}
P_{U_0}=\left(
\begin{array}{lcl}
s & p^{(2k+1)}(z) \\
 0 & t
\end{array}
\right)\text{ on $U_0$,}\quad
P_{U_{\infty}}=\left(
\begin{array}{lcl}
s & w^{2k+1}p^{(2k+1)}(1/w) \\
 0 & t
\end{array}
\right)\text{ on $U_{\infty}$}
\end{equation*}
where $s, t\in \bC$ and $p^{(2k+1)}(z)$ is a polynomial in $z$ of degree at most $2k+1$.
If $st\neq 0$, then $P\in \Aut(E)$.
We take $P\in \Aut(E)$.
Then we have
\begin{align}
 &\text{the $(1,1)$-entry of $P^{-1}_{U_0} A_z P_{U_0}$}  = \frac{t f_{11}^{(n-2)}(z) -p^{(2k+1)}(z)f_{21}^{(n-2k-3)}}{t}\\
 &\text{the $(2,1)$-entry of $P^{-1}_{U_0} A_z P_{U_0}$} = \frac{sf_{21}^{(n-2k -3)}(z)}{t}.
\end{align}

We consider simple descriptions of Higgs fields by the automorphisms of $E$. 
First, we consider the $(2,1)$-entry.
Let $\{ [s_1: 1],\ldots,[s_i:1], [ 1:q_{i+1}],\ldots,[1:q_{n-2k-3}] \}$ be the zeros of $f_{21}^{(n-2k-3)}(z)$ on $\mathbb{P}^1$
where $0\le i \le n-2k -3$.
By the automorphisms of $E$, we can put
\begin{equation*}
 f_{21}^{(n-2k-3)}(z) := (s_1 z-1)\cdots (s_i z-1) (z-q_{i+1}) \cdots (z-q_{n-2 k-3} ).
\end{equation*}
Second, we consider the $(1,1)$-entry.
We assume that the coefficient of $z^l$ in the polynomial $f_{21}^{(n-2k-3)}(z)$ is nonzero for some $l$ ($0\le l\le n-2k-3$).
By the automorphisms of $E$, we can put
\begin{equation}\label{normal form}
 f_{11}^{(n-2)}(z) := a_{n-2} z^{n-2}+\cdots +a_{l+2k+2} z^{l+2k+2}   + a_{l-1}z^{l-1}+\ldots+a_0.
\end{equation}
In particular, if $i=0$, that is, $ f_{21}^{(n-2k-3)}(z) := ( z-q_1) \cdots (z-q_{n-2 k-3} )$, 
then the coefficient of $z^{n-2k-3}$ in the polynomial $f_{21}^{(n-2k-3)}(z)$ is nonzero.
In this case, we can put
\begin{equation}\label{normal form n-2k-4}
 f_{11}^{(n-2)}(z) :=   a_{n-2k-4}z^{n-2k-4}+\ldots+a_0.
\end{equation}

\subsection{Apparent singularities and the dual parameters}\label{Subsection App}
We recall \textit{the apparent singularities and the dual parameters} introduced by Saito-Szabo \cite{SS}.
Let $(E,\Phi,\varphi) \in M_{{\it H}}$.
If $E \cong \cO\oplus \cO(-1)$,
then the apparent singularities and the dual parameters coincide with the geometric Darboux coordinates due to Oblezin \cite[Section 3]{Obl},
which gives a geometric interpretation of the Sklyanin formulas from \cite{Skl}.
We fix a section $\sigma \in H^0(\bP^1, E)$.
For the section $\sigma$, we define the following composition
\begin{equation*}
\cO_{\bP^1} \xrightarrow{\ \sigma\ } E \xrightarrow{\ \Phi\ } E \otimes L \lra (E/\cO_{\bP^1}) \otimes L.
\end{equation*}
The composition $\cO_{\bP^1}\ra (E/\cO_{\bP^1}) \otimes L$ is injective.
Then we can define a subsheaf $F^0\subset E $ such that $\cO_{\bP^1} \ra (F^0/\cO_{\bP^1}) \otimes L$ is isomorphic.
By the isomorphism $F^0/\cO_{\bP^1} \cong L^{-1}$, we have $F^0 \cong \cO_{\bP^1} \oplus  L^{-1}$.
Therefore, we have the following exact sequence.
\begin{equation}\label{ES of App}
0 \lra \cO_{\bP^1} \oplus L^{-1} \lra E \lra T_A \lra 0
\end{equation}
where $T_A$ is a torsion sheaf.
By the Riemann-Roch theorem, we have that the torsion sheaf $T_A$ is length $n-3$.
The exact sequence (\ref{ES of App}) is called a Frobenius--Hecke sheaf
originally introduced by Drinfeld (see \cite{Drinfeld} and \cite[Section 3.3]{Obl}).

\begin{Def}
For $(E,\Phi , \varphi) \in M_H$ and a nonzero section $\sigma \in H^0(\bP^1,E)$, 
we call the support of $T_A$ {\it apparent singular points of a $\bnu$-$\mathfrak{sl}_2$-parabolic Higgs bundle with a cyclic vector $(E,\Phi , \varphi,[\sigma])$}.
\end{Def}

Next, we define \textit{dual parameters of} $(E,\Phi , \varphi,[\sigma])$.
Let $C_s$ be the spectral curve of $(E,\Phi,\varphi)$.
Let $G$ be a torsion free sheaf of rank 1 on $C_s$ corresponding to $(E,\Phi,\varphi)$, which satisfies $E =\pi_*G$.
Since $H^0 (C_s, G) \cong H^0(\mathbb{P}^1, E)$, for a section $\sigma \in H^0(\mathbb{P}^1, E)$, 
we have the short exact sequence 
\begin{equation*}
0 \lra \cO_{C_s} \xrightarrow{\ \sigma\ } G \lra T_B \lra 0
\end{equation*}
where $T_B$ is a torsion sheaf on $C_s$ of length $n-3$.
We take the direct image of the short exact sequence.
Since $\pi_*(\cO_{C_s}) = \cO_{\bP^1} \oplus L^{-1}$ and $\pi_*G=E$, we have
\begin{equation*}
0 \lra  \cO_{\bP^1} \oplus L^{-1} \xrightarrow{\ \pi_*\sigma\ } E \lra \pi_*(T_B) \lra 0.
\end{equation*}
We may show that this short exact sequence coincides with the short exact sequence (\ref{ES of App}).
In particular, we have $\pi_*(T_B) = T_A$, whose support is the apparent singularities of $(E,\Phi,\varphi,[\sigma])$.
Set $\mathrm{Supp}(T_B)=\{ \tilde{p}_1, \ldots, \tilde{p}_{n-3} \}$, where $\tilde{p}_{i}$ is a point on $C_s$.
Put $\tilde{p}_i = ( q_i, p_i )$ where $q_i = \pi(\tilde{p}_i)$, which is an apparent singularity, and $p_i \in L_{q_i}$.

\begin{Def}
For $(E,\Phi , \varphi) \in M_H$ and a nonzero section $\sigma \in H^0(\bP^1,E)$, 
we call $\{ p_1,\ldots,p_{n-3} \}$ {\it dual parameters of a $\bnu$-$\mathfrak{sl}_2$-parabolic Higgs bundle with a cyclic vector $(E,\Phi , \varphi,[\sigma])$}.
\end{Def}

We consider the apparent singularities and the dual parameters of $(E,\Phi,\varphi,[\sigma])$ where $E\cong \cO\oplus \cO(-1)$.
In this case, the Higgs field $\Phi$ is described as follows:
\begin{equation*}
\Phi = 
\begin{cases}
 A_z^0 \otimes \omega_z    & \text{on } U_0 \\
 R_0^{-1}(A_z^0 \otimes \omega_z )R_0  & \text{on } U_{\infty}
\end{cases}\quad \text{where }\ 
A_z^0 :=\left(
\begin{array}{clcl}
f_{11}^{(n-2)}(z) & f_{12}^{(n-1)}(z) \\
 f_{21}^{(n-3)}(z) & -f_{11}^{(n-2)}(z)
\end{array}
\right).
\end{equation*}
The apparent singularities are the zeros of $f_{21}^{(n-3)}(z)$ on $\bP^1$, denoted by $\{ q_1 ,\ldots,q_{n-3} \}$.
The dual parameters are $\{ p_1, \ldots, p_{n-3}\}$ where we put $p_i:= f_{11}^{(n-2)}(q_i)$.
Then, for $(\cO\oplus \cO(-1),\Phi,\varphi)$, we have 
\begin{equation*}
\{ (q_1,p_1),\ldots, (q_{n-3},p_{n-3}) \} \in \mathrm{Sym}^{n-3}(\cK_n).
\end{equation*}

Next, we consider the case $E\cong \cO(k)\oplus \cO(-k-1)$ where $k>0$.
In this case, the Higgs field $\Phi$ is described as follows:
\begin{equation*}
\Phi = 
\begin{cases}
 A_z^k \otimes \omega_z    & \text{on } U_0 \\
 R_0^{-1}(A_z^k \otimes \omega_z )R_0  & \text{on } U_{\infty}
\end{cases}\quad \text{where }\ 
A_z^k :=\left(
\begin{array}{clcl}
f_{11}^{(n-2)}(z) & f_{12}^{(n+2k-1)}(z) \\
 f_{21}^{(n-2k-3)}(z) & -f_{11}^{(n-2)}(z)
\end{array}
\right).
\end{equation*}
Let $\sigma \in H^0(\bP^1,E) \cong H^0(\bP^1, \cO(k))$ be a section of $E$,
and let $\{ q_1,\ldots, q_k \}\in \mathrm{Sym}^k(\bP^1)$ be the zeros of the section $\sigma$.
We denote by $\{ q_{2k+1} ,\ldots,q_{n-3} \}\in \mathrm{Sym}^{n-2k-3}(\bP^1)$ the zeros of $f_{21}^{(n-2k-3)}(z)$ on $\bP^1$.
The apparent singularities of $(E,\Phi, \varphi, [\sigma])$ are the following
\begin{equation*}
\{ 2q_1, \ldots,2q_k,q_{2k+1},\ldots, q_{n-3} \} \in \mathrm{Sym}^{n-3}(\bP^1).
\end{equation*}
We compute the dual parameters of $(E,\Phi, \varphi, [\sigma])$.
We take $\sigma \in H^0(C_s , G) \cong H^0(\bP^1,\pi_*G)$ corresponding to $\sigma \in H^0(\bP^1,E)$.
The section $\sigma \in H^0(C_s , G)$ has the following zeros 
\begin{equation*}
\{ (q_1, p_1),(q_1, -p_1),\ldots , (q_k,p_k),(q_k,-p_k) \}
\end{equation*}
where $\det ( p_i I- \Phi|_{q_i})=0$ for $i=1,\ldots,k$.
The dual parameters are $\{ p_1, -p_1, \ldots, p_k,-p_k, p_{2k+1}, \ldots, p_{n-3}\}$ where we put $p_i:= f_{11}^{(n-2)}(q_i)$ for $i=2k+1 ,\ldots, n-3$.
Then, for $(E,\Phi,\varphi,[\sigma])$, we have 
\begin{equation*}
\{ (q_1, p_1),(q_1, -p_1),\ldots , (q_k,p_k),(q_k,-p_k), (q_{2k+1},p_{2k+1}),\ldots, (q_{n-3},p_{n-3}) \} \in \mathrm{Sym}^{n-3}(\cK_n).
\end{equation*}

\subsection{Lower and upper modifications}
In this subsection, following \cite[Section 2]{Obl}, we describe \textit{the lower and the upper modifications}.
Let $E$ be an algebraic vector bundle on $\bP^1$ of rank $2$ and of degree $d$.
Fix a point $a \in \bP^1$.
Let $l \subset E|_a$ be a $1$-dimensional subspace. 
\begin{Def}
We call  
\begin{equation*}
(a,l)^{\text{low}}(E) := \{ s \in E \mid s(a) \in l \}, \quad (a,l)^{\text{up}}(E) := (a,l)^{\text{low}}(E) \otimes \cO(a)
\end{equation*}
the lower and the upper modifications of $E$, respectively.
\end{Def}
\textit{The lower and the upper modifications} provide the following exact sequences
\begin{equation*}
0\lra (a,l)^{\text{low}}(E) \lra E \lra E|_{a}/l \lra 0,
\end{equation*}
\begin{equation*}
0\lra E \lra (a,l)^{\text{up}}(E) \lra l \otimes \cO(a) \lra 0,
\end{equation*}
respectively.
In other words, we change our bundle rescaling the basis of sections in the neighborhood of a point $a$ as follows.
Given a local decomposition $V=l\oplus l'$ of $E\cong V \otimes \cO$,
we put the local basis $\{ s_1(z), s_2(z) \}$ with $l\otimes \cO \cong \langle s_1(z) \rangle$ and $l'\otimes \cO \cong \langle s_2(z) \rangle$.
Then the basis of the lower modification $(x,l)^{\text{low}}$ of the bundle is generated by the sections $\{ s_1(z),(z-x) s_{2}(z)\}$,
and of the upper one $(x,l)^{\text{up}}$ by $\{ (z-x)^{-1}s_1(z), s_{2}(z)\}$.
Consequently, in the punctured neighborhood, we may represent the action of the modifications by the following gluing matrices
\begin{equation*}
(a,l)^{\text{low}}=\left(
\begin{array}{ll}
1 & 0 \\
 0 & (z-a)
\end{array}
\right),\quad
(a,l)^{\text{up}}=\left(
\begin{array}{ll}
(z-a)^{-1}  & 0 \\
 0 & 1
\end{array}
\right).
\end{equation*}

\section{Geometric description of the moduli spaces}
Suppose that $\bnu$ satisfies the condition (\ref{generic condition 2}) and $\nu_1\cdots\nu_n\neq0$.
We put 
\begin{equation*}
\begin{aligned}
(t_1,t_2,t_{n})&:=(0,1,\infty), \\
(\nu^{\pm}_{1},\ldots,\nu^{\pm}_{n})&:=(\pm\nu_1,\pm\nu_2,\ldots,\pm\nu_{n}), \text{ and}\\
\hat{\nu}_i &:=  \nu_i (t_i-t_1)\cdots(t_i-t_{i-1})(t_i-t_{i+1})\cdots (t_i-t_{n-1}) \text{ for $i=1,\ldots,n-1$}.
\end{aligned}
\end{equation*}
First, we consider the apparent singularities and the dual parameters of members of $M_{{\it H}}^0$ for $n\ge4$.
Then we have Theorem \ref{Prop M^0 to Hilb 0} (Theorem \ref{Prop M^0 to Hilb}). 
Second, we assume $n=5$.
We consider the apparent singularities and the dual parameters of members of $\widehat{M}_{{\it H}}$.
Then we have the first assertion of Theorem \ref{main theorem} (Theorem \ref{main theorem +}).

\subsection{Geometric description of $M_{{\it H}}^0$ for $n\ge4$}\label{O+O(-1)}

Let $(\cO\oplus\cO(-1),\Phi, \varphi) \in M_{{\it H}}^0$, and $\cK_n'$ be the Zariski open set of the blowing-up of Hirzebruch surface of degree $n-2$ defined in \ref{Hirzebruch blow up},
and $\cK_n$ be the contraction $\cK'_n\ra \cK_n$. 
Since $\dim H^0(\bP^1,\cO\oplus\cO(-1))=1$, sections are determined uniquely up to constant.
Then the apparent singularities and the dual parameters are determined by $(E,\Phi,\varphi)$.
Let
$\{ (q_1,p_1), \ldots,(q_{n-3},p_{n-3})  \} \in \mathrm{Sym}^{n-3}(\cK_n)$
be the apparent singularities and the dual parameters of $(E,\Phi,\varphi)$.
We consider the map
\begin{equation}\label{M to Sym}
\begin{aligned}
M^0_{{\it H}} &\lra \mathrm{Sym}^{n-3}(\cK_n)\\
(E,\Phi,\varphi) &\longmapsto  \{ (q_1,p_1), \ldots,(q_{n-3},p_{n-3})  \},
\end{aligned}
\end{equation}
which is essentially constructed in \cite[Section 3]{Obl}.
Since this map is not injective, 
we consider the composite of the Hilbert-Chow morphism and the blowing-up
\begin{equation*}
\Hilb^{n-3}(\cK'_n) \lra \mathrm{Sym}^{n-3}(\cK'_n) \lra \mathrm{Sym}^{n-3}(\cK_n).
\end{equation*}
Then we have the following 
\begin{Thm}\label{Prop M^0 to Hilb}
\textit{
The map} (\ref{M to Sym}) \textit{is extended to 
\begin{equation}
M_{{\it H}}^0 \lra \Hilb^{n-3}(\cK'_n).
\end{equation}
The map is injective.
Moreover, we can give an explicit description of the universal family $(\tilde{E}^{(0)}, \tilde{\Phi}^{(0)}) \ra M_{{\it H}}^0\times \mathbb{P}^1$.
}
\end{Thm}
The image of the map $M_{{\it H}}^0 \ra \Hilb^{n-3}(\cK'_n)$ is described as follows.
The image of the map $M_{{\it H}}^0 \ra \mathrm{Sym}^{n-3}(\cK'_n)$ is 
\begin{equation*}
\{ x=\{ n_1 \tilde{p}_1, \ldots,n_r \tilde{p}_r \} \in  \mathrm{Sym}^{n-3}(\cK_{n}') \mid \pi(\tilde{p}_i) \neq \pi(\tilde{p}_j) \text{ for } i\neq j\}
\end{equation*}
where $n_j$ are integers such that $n_1+\cdots+n_{r}=n-3$, $n_j\ge 1$ ($j=1,\ldots,r$) and $\pi$ is the projection $\cK_{n}' \ra \mathbb{P}^1$.
For $x=\{ n_1 \tilde{p}_1, \ldots,n_r \tilde{p}_r \}$, let $M_{H,x}^0$ be the fiber of $x$ under $M_{{\it H}}^0 \ra \mathrm{Sym}^{n-3}(\cK'_n)$.
Then $M_{{\it H}}^0 \cong \mathbb{C}^{n_1-1}\times \cdots \times \mathbb{C}^{n_r-1}$ 
where $\mathbb{C}^{n_1-1}\times \cdots \times \mathbb{C}^{n_r-1}$ is an affine open set of the fiber of $x$ under the Hilbert-Chow morphism.

The explicit description of the family is as follows.
For simplicity, we assume that $\{ q_1,\ldots,q_{n-3} \}\subset  \mathbb{P}^1 \setminus \{ \infty \}$.
If the apparent singularities are distinct, then the explicit description is the following
\begin{equation}
\begin{cases}
 A_z^0 \otimes \omega_z    & \text{on } U_0 \\
 R_0^{-1}(A_z^0 \otimes \omega_z )R_0  & \text{on } U_{\infty}
\end{cases},\quad
A_z^0= 
\begin{pmatrix}
a_{n-4}z^{n-4} + \cdots+ a_0 & f_{12}^{(n-1)}(z) \\
(z-q_1)\cdots (z-q_{n-3}) & -(a_{n-4}z^{n-4} + \cdots+ a_0)
\end{pmatrix}
\end{equation}
where
\begin{equation}\label{family coefficient 00}
a_i = \sum_{k=1}^{n-3} \sum_{j=1}^{n-3} u_{i+1, j} l_{jk} p_{k}, \quad R_0=
\begin{pmatrix}
1 & 0\\
0 & \frac{1}{z}.
\end{pmatrix}
\end{equation}
Here $(q_1,p_1),\ldots, (q_{n-3},p_{n-3})$ are pairs of the apparent singularities and their duals,
the element $l_{ij}$ is defined by the relations $l_{ij}=0$ for $i<j$, $l_{11}=1$ and $l_{ij}= \prod_{k=1, k\neq j}^i 1/(q_j-q_k)$ otherwise,
and the element $u_{ij}$ is defined by the relations $u_{ii}=1$, $u_{i1}=0$ and $u_{ij}= u_{i-1,j-1} - u_{i,j-1}q_{j-1}$ otherwise,
that is, the matrix $(u_{ij} l_{jk})_{ik}$ is the inverse of the Vandermonde matrix.
Here we set $u_{0j}=0$.
We omit the description of $f_{12}^{(n-1)}(z)$, which is a polynomial in $z$ of degree $n-1$, since the description is lengthened and is not necessary below.
Next, we consider the case where the apparent singularities and their duals have multiplicities.
Let $n_j$ ($j=1,\ldots,r$) be integers such that $n_1+\cdots+n_{r}=n-3$, $n_j\ge 1$.
For each $j$ ($j=1,\ldots,r$), we assume that $q_i=x_j \in \mathbb{P}^1\setminus \{ \infty \}$ and $p_i =y_j \in \mathbb{C}$ ($i=\sum_{k=1}^{j-1} n_k+1, \cdots, \sum_{k=1}^{j} n_k$), that is, 
the pair of the apparent singularities and their duals is $\{ n_1 (x_1, y_1), \ldots,n_r (x_r, y_r) \}$.
In this case, for each $j$ ($j=1,\ldots,r$) we substitute 
\begin{equation*}
p_i=y_j +(q_i-x_j)(\lambda_0^{(j)} +\lambda_1^{(j)} (q_i-x_j) + \cdots + \lambda_{n_j-2}^{(j)} (q_i-x_j)^{n_j-2}), \quad i=\sum_{k=1}^{j-1} n_k+1, \cdots, \sum_{k=1}^{j} n_k,
\end{equation*}
for the coefficient (\ref{family coefficient 00}).
Here $\lambda_0^{(j)},\ldots, \lambda_{n_j-2}^{(j)}$ are coordinates of the affine open set $\mathbb{C}^{n_j-1}$ of the fiber of the Hilbert-Chow morphism.
Then we can define the coefficients $a_i$ for this case.

\begin{proof}[Proof of Theorem \ref{Prop M^0 to Hilb}]
Put
\begin{equation*} 
M^{00}_{{\it H}} := \left\{ (E,\Phi,\varphi) \in M^{0}_{{\it H}}\ \middle|  
\begin{array}{l}
\{ q_1,\ldots,q_{n-3} \}: \text{the apparent singularities of } (E,\Phi,\varphi)\\
\ q_i\neq q_j \ (i \neq j)
\end{array}
  \right\},
\end{equation*}
and
\begin{equation*} 
M^{000}_{{\it H}} := \left\{ (E,\Phi,\varphi) \in M^{0}_{{\it H}}\ \middle|
\begin{array}{l}
\{ q_1,\ldots,q_{n-3} \}: \text{the apparent singularities of } (E,\Phi,\varphi)\\
\ q_i\neq q_j \ (i \neq j) \text{ and } t_k \notin \{ q_1,\ldots,q_{n-3} \} \text{ for any $k$}.
\end{array}
  \right\}.
\end{equation*}

\underline{{\bf Step 1}}. In this step, we show that the restriction $M_{{\it H}}^{000} \rightarrow \mathrm{Sym}^{n-3}(\cK_n)$ is injective. 
The image of this restriction is the following
\begin{equation*}
\mathrm{Image}(M^{000}_{{\it H}}):= \left\{ \{ (q_1,p_1) ,\ldots, (q_{n-3},p_{n-3} ) \}\ \middle|
\begin{array}{l}
q_i \neq q_j\ (i\neq j), \\
q_i \notin \{ t_1,\ldots,t_n \}\ (i=1,\ldots,n-3)  
\end{array}
\right\} \subset \mathrm{Sym}^{n-3}(\cK_n).
\end{equation*}
Let 
\begin{equation}\label{element of image 1}
\{ ( [s_1:1],u_1),\ldots,([s_i:1] ,u_i ),( [1:q_{i+1}],p_{i+1} ),\ldots,([1:q_{n-3}], q_{n-3} )  \} 
\end{equation}
be an element of $\mathrm{Image}(M^{000}_{{\it H}})$ where $0\le i \le n-3$.
We show that the entries $f_{21}^{(n-3)}(z)$, $f_{11}^{(n-2)}(z)$, and $f_{12}^{(n-1)}(z)$ of the description (\ref{normal form of O+O(-1)}) are determined by 
the element (\ref{element of image 1}) up to automorphisms of $\cO\oplus \cO(-1)$.
First, we consider the entry $f_{21}^{(n-3)}(z)$.
By the definition of the apparent singularities and an automorphism of $\cO\oplus\cO(-1)$, we can put
\begin{equation*}
f_{21}^{(n-3)}(z) =  (s_1 z-1)\cdots (s_i z-1) ( z-q_{i+1}) \cdots ( z-q_{n-3}).
\end{equation*}

Second, we consider the entry $f_{11}^{(n-2)}(z)$.
Since $s_1\cdots s_i\neq 0$, the coefficient $z^{n-3}$ in $f_{21}^{(n-3)}(z)$ is nonzero. 
Then, by the automorphism of $\cO\oplus\cO(-1)$, we can put 
$f_{11}^{(n-2)}(z)=a_{n-4}z^{n-4}+\cdots  + a_{0}$
as in the description (\ref{normal form n-2k-4}).
By the definition of the dual parameters, we have that
$u_j= s_j^{n-2}f_{11}^{(n-2)}(1/s_j)$ for $1\le j \le i$ and $p_i= f_{11}^{(n-2)}(q_j)$ for $i+1\le j \le n-3$.
Then we have the following system
\begin{equation*}
\begin{pmatrix}
u_1\\
\vdots\\
u_i\\
p_{i+1}\\
\vdots\\
p_{n-3}\\
\end{pmatrix}
=
\begin{pmatrix}
1 & \cdots &s_{1}^{n-3} & s_{1}^{n-4}\\
  \vdots  &       &  \vdots     & \vdots \\
1 & \cdots &s_{i}^{n-3} & s_{i}^{n-4}\\
q_{i+1}^{n-4} & \cdots &q_{i+1} & 1\\
  \vdots  &       &  \vdots     & \vdots \\
q_{n-3}^{n-4} & \cdots  &q_{n-3} & 1
\end{pmatrix}
\begin{pmatrix}
a_{n-4}\\
\vdots\\
a_1\\
a_{0}
\end{pmatrix}.
\end{equation*}
We can determine the coefficients $a_{n-4}, \ldots a_0$ by an element (\ref{element of image 1}).

Third, we consider the entry $f_{12}^{(n-1)}(z)$.
We put
$f_{12}^{(n-1)}(z):=b_{n-1}z^{n-1}+b_{n-2}z^{n-2}+\cdots + b_0$.
We solve the equations
\begin{equation}\label{condition local exp}
\det\left( \res_{t_i} \Phi \right) = - \nu_i^2,\quad i=1,\ldots,n.
\end{equation}
Then we have 
\begin{equation*}
f_{12}^{(n-1)}(t_i)= \frac{-f_{11}^{(n-2)}(t_i)^2+\hat{\nu}_i^2}{f_{21}^{(n-3)} (t_i)} \ (\text{for }i=1,\ldots,n-1), \text{ and }
b_0=\frac{-(w^{n-2} f_{11}^{(n-2)}(1/w))^2|_{w=0}+\nu_{n}^2}{(w^{n-3} f_{21}^{(n-3)} (1/w) )|_{w=0}}.
\end{equation*}
Since $f_{11}^{(n-2)}$ and $f_{21}^{(n-3)}$ are determined by the element (\ref{element of image 1}),
we can determine the coefficients $b_0,\ldots,b_{n-1}$.
As the result, we obtain that the map $M_{{\it H}}^{000} \rightarrow \mathrm{Sym}^{n-3}(\cK_n)$ is injective.

\underline{{\bf Step 2}}. In this step, we extend the map $M_{{\it H}}^{000} \rightarrow \mathrm{Sym}^{n-3}(\cK_n)$ to
$M_{{\it H}}^{00} \rightarrow \mathrm{Sym}^{n-3}(\cK_n')$, and  
we show that the extended map $M_{{\it H}}^{00} \rightarrow \mathrm{Sym}^{n-3}(\cK_n')$ is injective. 
For the element (\ref{element of image 1}), we can put
\begin{equation*}
f_{11}^{(n-2)}(z):= p_j+  (z-q_j)\tilde{f}_{11}(z), \text{ and }f_{21}^{(n-3)}(z):= (s_1 z-1)\cdots (s_i z-1) ( z-q_{i+1}) \cdots ( z-q_{n-3})
\end{equation*}
where $\tilde{f}_{11}(z)$ is a polynomial of degree at most $n-3$ in $z$.
The polynomial $\tilde{f}_{11}(z)$ is determined by the element (\ref{element of image 1}) up to automorphisms of $\cO\oplus \cO(-1)$.
By the condition (\ref{condition local exp}), we have
\begin{equation}\label{f12ti}
\begin{aligned}
f_{12}^{(n-1)}(t_i) &=\frac{-(p_j+  (t_i-q_j)\tilde{f}_{11}(t_i))^2+  \hat{\nu}_i^2}{(t_i-q_1)\cdots (t_i-q_{n-3})} \\
&=\frac{-(t_i-q_j)(\tilde{f}_{11}(t_i)^2 + 2 p_j \tilde{f}_{11}(t_i)) -p_j^2+  \hat{\nu}_i^2 }{(t_i-q_1)\cdots (t_i-q_{n-3})} .
\end{aligned}
\end{equation}
We consider the blowing-up $\cK'_n\ra \cK_n$.
Let $\epsilon \in \{ +, -\}$.
We define the blowing-up parameters $v_j^{i,\epsilon}$ at $(t_i, \epsilon \hat{\nu}_i) \in \cK_n$ as $p_j - \epsilon \hat\nu_i = v_j^{i,\epsilon}(q_j-t_i)$.
We substitute $p_j = v_j^{i,\epsilon}(q_j-t_i) + \epsilon \hat{\nu}_i$ for the formula (\ref{f12ti}).
Then we have
\begin{equation}
\begin{aligned}
f_{12}^{(n-1)}(t_i) &=\frac{-(t_i-q_j)(\tilde{f}_{11}(t_i)^2 + 2 p_j \tilde{f}_{11}(t_i)) -(v_j^{i,\epsilon}(q_j-t_i))^2- \epsilon 2 \hat{\nu}_i v_j^{i,\epsilon}(q_j-t_i)}{(t_i-q_1)\cdots (t_i-q_{n-3})} \\
&=\frac{-(\tilde{f}_{11}(t_i)^2 + 2 p_j \tilde{f}_{11}(t_i)) -(v_j^{i,\epsilon})^2(q_j-t_i) + \epsilon 2 \hat{\nu}_i v_j^{i,\epsilon}}{(t_i-q_1)\cdots (t_i-q_{j-1})(t_i-q_{j+1}) \cdots (t_i-q_{n-3})}.
\end{aligned}
\end{equation}
We consider the behavior of $f_{12}^{(n-1)}(t_i)$ as $q_j \ra t_i$.
The limit $\lim_{q_j\ra t_i} f_{12}^{(n-1)}(t_i)$ is convergence, 
and the convergence value is determined by the apparent singularities, the dual parameters, and the blowing-up parameters $v_j^{i,\epsilon}$.
By the same argument as in Step 1, we can determine the all coefficients of $f_{12}^{(n-1)}(z)$.
Then we obtain the map $M_{{\it H}}^{00} \rightarrow \mathrm{Sym}^{n-3}(\cK_n')$, and this map is injective.

\underline{{\bf Step 3}}. In this step, we extend the map $M_{{\it H}}^{00} \rightarrow \mathrm{Sym}^{n-3}(\cK_n')$ to
$M_{{\it H}}^{0} \rightarrow \mathrm{Hilb}^{n-3}(\cK_n')$, and  
we show that the extended map $M_{{\it H}}^{0} \rightarrow \mathrm{Hilb}^{n-3}(\cK_n')$ is injective.
Let 
\begin{equation}\label{element of image 2}
\{ ( [s_1:1],u_1),\ldots,([s_{n_0}:1] ,u_{n_0} ),( [1:q_{1}^1],p_{1}^1 ),\ldots,( [1:q_{n_1}^{1}],p_{n_1}^1 ),\ldots ,( [1:q_{1}^r],p_{1}^r ),\ldots,( [1:q_{n_r}^{r}],p_{n_r}^r )  \} 
\end{equation}
be an element of $\mathrm{Image}(M^{000}_{{\it H}})$ where $n_0+\cdots+n_{r}=n-3$, $n_j\ge 1$ ($j=1,\ldots,r$).
By the element (\ref{element of image 2}), we can determine the entries $f_{21}^{(n-3)}(z)$, $f_{11}^{(n-2)}(z)$, and $f_{12}^{(n-1)}(z)$ of the description (\ref{normal form of O+O(-1)}) (Step 1).
Let $\tilde{x}$ be a point of $\mathrm{Sym}^{n-3}(\cK_{n}')$, denoted by
\begin{equation*}
\tilde{x}= \{ n_0 ([0:1], y_0,a_0^{i,\pm}), n_1 ([1:x_1], y_1,a_1^{i,\pm}), \ldots,n_r ([1:x_{r}],y_{r},a_r^{i,\pm}) \} \in  \mathrm{Sym}^{n-3}(\cK_{n}')
\end{equation*}
where $y_l+ \epsilon \hat{\nu}_i = a_l^{i,\epsilon} (x_l -t_i)$ for $l=0,\ldots,r$ and $\epsilon\in \{ +,-\}$.
We consider the behavior of the entries as $s_j\ra 0$ ($j=1,\ldots,n_0$), $q_k^l \ra x_l$, and $p_k^l \ra y_l$ ($k=1,\ldots, n_l$ and $l=1,\ldots,r$) where $x_{l_1}\neq x_{l_2}$ ($l_1\neq l_2$).
Let $\cI$ be an ideal contained in the fiber of $\tilde{x}$ by the Hilbert-Chow morphism $\Hilb^{n-3}(\cK'_n) \ra \mathrm{Sym}^{n-3}(\cK_{n}')$.
We show that the entries $f_{11}$, $f_{12}$, $f_{21}$ are determined by the ideal $\cI$ up to automorphisms.

On a neighborhood of $\cI$, the Hilbert scheme of points $\Hilb^{n-3}(\cK'_n)$ is isomorphic to 
$\Hilb^{n_0}(\cK'_n) \times \cdots \times \Hilb^{n_r}(\cK'_n)$.
We denote by $(\cI_{x_0,y_0}^{n_0}, \ldots,\cI_{x_r,y_r}^{n_r})$ the image of $\cI$.
We put $\cI_{x_l,y_l}^{n_l}=((q -x_l)^{n_l}, f_{x_l,y_l}(q,p) )$
where
$f_{x_l,y_l}(q,p) := (p-y_l) -(q-x_l)(\lambda_0^{(l)} +\lambda_1^{(l)} (q-x_j) + \cdots + \lambda_{n_l-2}^{(l)} (q-x_l)^{n_l-2})$
as in the proof of \cite[Theorem 1.13]{Naka}.

We consider the entry $f_{11}^{(n-2)}(z)$.
For a neighborhood of the ideal $\cI$, we can assume that the coefficient of $z^{n-n_0-3}$ in $f_{21}^{(n-3)}(z)$ is nonzero.
Then we can normalize $f_{11}^{(n-2)}(z)$ as
\begin{equation}\label{f11 step3} 
f_{11}^{(n-2)}(z) := a_{n-2} z^{n-2}+\cdots +a_{n-n_0-1} z^{n-n_0-1}   + a_{n-n_0-4}z^{n-n_0-4}+\ldots+a_0
\end{equation}
by automorphisms of $\cO\oplus \cO(-1)$.
By the definition of dual parameters, we have the following system
\begin{equation}\label{Vandermonde}
\begin{pmatrix}
u_1\\
\vdots\\
u_{n_0}\\
p^1_{1}\\
\vdots\\
p^r_{n_r}
\end{pmatrix}
=
\begin{pmatrix}
 1 & s_1 & \cdots & s_1^{n_0-1} & s_1^{n_0+2} &\cdots & s_1^{n-2}      \\
  \vdots  &  \vdots &    &  \vdots     & \vdots& & \vdots  \\
  1 & s_{n_0} & \cdots & s_{n_0}^{n_0-1} & s_{n_0}^{n_0+2} &\cdots & s_{n_0}^{n-2}      \\
(q^1_{1})^{n_1-2}& (q^1_{1})^{n-3} & \cdots &(q^1_{1})^{n-n_0-1} & (q^1_{1})^{n-n_0-4} & \cdots & 1\\
  \vdots  &  \vdots &    &  \vdots     & \vdots& & \vdots \\
(q^r_{n_r})^{n-2}& (q^r_{n_r})^{n-3} & \cdots &(q^r_{n_r})^{n-n_0-1} & (q^r_{n_r})^{n-n_0-4} & \cdots & 1
\end{pmatrix}
\begin{pmatrix}
a_{n-2}\\
\vdots\\
a_{n-n_0-1}\\
a_{n-n_0-4}\\
\vdots\\
a_{0}
\end{pmatrix}.
\end{equation}
We substitute 
\begin{equation*}
u_k=y_0 +s_k(\lambda_0^{(k)} +\lambda_1^{(k)} s_k + \cdots + \lambda_{n_0-2}^{(k)} s_k^{n_0-2}),\quad k=1 ,\ldots, n_0,
\end{equation*}
and
\begin{equation*}
p_k^l=y_l +(q_k^l-x_l)(\lambda_0^{(l)} +\lambda_1^{(l)} (q_k^l-x_l) + \cdots + \lambda_{n_l-2}^{(l)} (q_k^l-x_l)^{n_l-2}),\quad  k=1 ,\ldots, n_l,
\end{equation*}
where $j=1,\ldots,r$, in the system (\ref{Vandermonde}).
Then we can show that the coefficients of the polynomial (\ref{f11 step3}) are defined when $s_k=0$ and $q_k^l=x_l$.
Moreover, these coefficients are determined by the apparent singularities, the dual parameters,
and the parameters $\lambda_0^{(i)} ,\lambda_1^{(i)} , \ldots , \lambda_{n_i-2}^{(i)}$.

We consider the entry $f_{12}^{(n-1)}(z)$.
If $x_l \notin \{ t_1,\ldots,t_n\}$, then we have the value $f_{12}^{(n-1)}(t_i) \in \mathbb{C}$ as $q_k^l \ra x_l$ for any $k=1,\ldots,n_l$.
By the values $f_{12}^{(n-1)}(t_1),\ldots, f_{12}^{(n-1)}(t_n)$, we can determine all coefficients of $f_{12}^{(n-1)}(z)$.
Next, we consider the behavior of the value $f_{12}^{(n-1)}(t_i)$ as $x_l \ra t_i$.
Let $\epsilon \in \{ +,-\}$.
We consider the ideal $\cI_{x_l,a_l^{i,\epsilon}}^{n_l}=((q -t_i)^{n_l}, f_{x_l,a_l^{i,\epsilon}}(q,v^{i,\epsilon}) )$ where
\begin{equation*}
f_{x_l,a_l^{i,\epsilon}}(q,v^{i,\epsilon}) := (v^{l,\epsilon}-a^{i,\epsilon}) - (q - x_l)( \lambda_0^{(l)}  + \lambda_1^{(l)} (q-x_l) + \cdots  + \lambda_{n_l-2}^{(l)} (q-x_l)^{n_l-2}).
\end{equation*}
For the ideal $\cI_{x_l,a_l^{i,\epsilon}}^{n_l}$, we can describe the dual parameter $p_k^l$ ($k=1,\ldots,n_l$) as follows:
\begin{equation*}
\begin{aligned}
p_k^l &=\epsilon \hat{\nu}_i+ v_k^{i,\epsilon}(q_k^l-t_i)\\ 
&=\epsilon \hat{\nu}_i+ \left( a_l^{i,\epsilon} + (q_k^l - x_l)( \lambda_0^{(l)}  + \lambda_1^{(l)} (q_k^l-x_l) + \cdots  + \lambda_{n_l-2}^{(l)} (q_k^l-x_l)^{n_l-2}) \right) (q_k^l-t_i);\\
&=\epsilon \hat{\nu}_i+a_l^{i,\epsilon}(x_l-t_i )+ \hat{\lambda}_1 (q_k^l- x_l)+\hat{\lambda}_2 (q_k^l- x_l)^2 +\cdots+\hat{\lambda}_{n_l-1} (q_k^l- x_l)^{n_l-1} 
+ \lambda_{n_l-2}^{(l)} (q_k^l- x_l)^{n_l}
\end{aligned}
\end{equation*}
where $\hat{\lambda}_j:= \lambda_{j-2}^{(l)}+\lambda_{j-1}^{(l)}(x_l-t_i) $ for $j=1,\ldots, n_i-1$.
Here we put $\lambda_{-1}^{(l)}:= a_l^{i,\epsilon}$.
For the ideal $\cI_{x_l,a_l^{i,\epsilon}}^{n_l}$, we put
\begin{equation}\label{f11f21}
\begin{aligned}
f_{11}^{(n-2)}(z) &:=\epsilon \hat{\nu}_i+a_l^{i,\epsilon}(x_l-t_i )+ \hat{\lambda}_1 (z- x_l)+\hat{\lambda}_2 (z- x_l)^2
+\cdots+\hat{\lambda}_{n_l-1} (z- x_l)^{n_l-1} +(z-x_l)^{n_l} \tilde{f}_{11}(z); \\
f_{21}^{(n-3)}(z) &:= (z-x_l)^{n_l} (z-q_{n_i+1}) \cdots (z-q_{n-3}).
\end{aligned}
\end{equation}
We substitute $z=t_i$ for $f_{11}^{(n-2)}(z)$.
Then we have
$f_{11}^{(n-2)}(t_i) =\epsilon \hat{\nu}_i +(t_i-x_l)^{n_l} (\tilde{f}_{11}(t_i) - \lambda_{n_l-2}^{(l)})$.
We consider the value of $f_{12}^{(n-1)}(t_i)$ as follows:
\begin{equation*}
\begin{aligned}
f_{12}^{(n-1)}(t_i)&= \frac{-f_{11}^{(n-2)}(t_i)^2+\hat{\nu}_i^2 }{(t_i-x_l)^{n_l} (t_i-q_{n_l+1}) \cdots (t_i-q_{n-3})}= \frac{\alpha+ \beta(t_i-x_l)^{n_l}}{ (t_i-q_{n_l+1}) \cdots (t_i-q_{n-3})}
\end{aligned}
\end{equation*}
where $\alpha:=\epsilon 2\hat{\nu}_i  (\tilde{f}_{11}(q_j)- \lambda_{n_l-2}^{(l)} ) \text{ and } \beta :=\ (\tilde{f}_{11}(q_j)- \lambda_{n_l-2}^{(l)}  )^2$.
Then we have the finite value $f_{12}^{(n-1)}(t_i)$ as $x_l \ra t_i$.
By the values $f_{12}^{(n-1)}(t_1),\ldots, f_{12}^{(n-1)}(t_n)$ as $x_l \ra t_i$, we can determine all coefficients of $f_{12}^{(n-1)}(z)$.

We obtain an extended map $M_{{\it H}}^{0} \ra \mathrm{Hilb}^{n-3}(\cK_n')$
by the assignment of $\cI_{x_l,a_l^{i,\epsilon}}^{n_l}=((q -t_i)^{n_l}, f_{x_l,a_l^{i,\epsilon}}(q,v^{i,\epsilon}) )$ to the matrix 
$\begin{pmatrix}
f^{(n-2)}_{11} & f^{(n-1)}_{12} \\
f^{(n-3)}_{21} & f^{(n-2)}_{22}
\end{pmatrix}$.
This procedure is confirmed by the implication that the existence of $\lim_{x_l\ra t_i} f_{12}^{(n-1)}(t_i)$ determines the data $f_{x_l,a_l^{i,\epsilon}}(q, v^{i,\epsilon})$ 
from the argument above.
The extended map is injective.
\end{proof}

\subsection{Geometric description of $\widehat{M}_{{\it H}}$ for $n=5$}\label{O(1)+O(-2)}
Suppose that $n=5$.
By the apparent singularities and the dual parameters of $(E,\Phi,\varphi,[\sigma]) \in \widehat{M}_{{\it H}}$,
we have the following map
\begin{equation}\label{Msys to Sym}
\begin{aligned}
\widehat{M}_{{\it H}} &\lra \mathrm{Sym}^2(\cK_5)\\
(E,\Phi,\varphi,[\sigma]) &\longmapsto \{ (q_1,p_1),(q_2,p_2)\}
\end{aligned}
\end{equation}
where $\{q_1, q_2\}$ are apparent singularities and $\{ p_1, p_2\}$ are their dual parameter where $p_i$ corresponds to $q_i$ for $i=1,2$.
There exists a stratification 
$\widehat{M}_{{\it H}} = \widehat{M}^{0}_{{\it H}} \cup \widehat{M}^{1}_{{\it H}}$
where $\widehat{M}^{i}_{{\it H}}$ is the locus such that $(E,\Phi,\varphi,[\sigma]) \in \widehat{M}^{i}_{{\it H}}$ satisfies $E \cong \cO(i)\oplus\cO(-i-1)$.
The image of $\widehat{M}^{1}_{{\it H}}$ is the following
\begin{equation*}
\mathrm{Image}(\widehat{M}^{1}_{{\it H}}):= \{ \{ (q_1,p_1) , (q_2,p_2 ) \} \mid q_1 = q_2,\ p_1=-p_2  \} \subset \mathrm{Sym}^{n-3}(\cK_5).
\end{equation*}

We take a blowing-up of $\mathrm{Hilb}^2(\cK'_5)$ as follows.
Let $Z \subset \mathrm{Sym}^2(\cK'_5)$ be the proper pre-image of $\{ (q_1,p_1),(q_1,-p_1) \} \subset \mathrm{Sym}^2 (\cK_5)$
for $\mathrm{Sym}^2(\cK'_5) \ra \mathrm{Sym}^2(\cK_5)$.
Let $\widetilde{Z} \subset \Hilb^2(\cK'_5)$ be the proper pre-image of $Z$
for $\Hilb^2(\cK'_5) \ra \mathrm{Sym}^2(\cK'_5)$.
We denote by 
\begin{equation}\label{blowing up along Z}
\widetilde{\mathrm{Hilb}}^2(\cK'_5) \lra \mathrm{Hilb}^2(\cK'_5)
\end{equation}
the blowing-up along $\widetilde{Z}$.
Then we have the following 
\begin{Thm}\label{main theorem +}
\textit{
The map} (\ref{Msys to Sym}) \textit{is extended to 
\begin{equation}\label{equ of main theorem +}
\widehat{M}_{{\it H}} \lra \widetilde{\mathrm{Hilb}}^2(\cK'_5).
\end{equation}
The map is injective.
}
\end{Thm}

\begin{proof}
Since $\widehat{M}^{0}_{{\it H}}=M^{0}_{{\it H}}$, 
we have the injective map $\widehat{M}^{0}_{{\it H}} \ra \mathrm{Hilb}^2(\cK'_5)$ by Theorem \ref{Prop M^0 to Hilb}. 
Set
\begin{equation*}
\widehat{M}^{10}_{{\it H}}:= \left\{ (E,\Phi,\varphi,[\sigma]) \in \widehat{M}^{1}_{{\it H}}\ \middle|
\begin{array}{l}
\text{the apparent singularities $\{ q_1 ,q_2 \}$ of } (E,\Phi,\varphi,[\sigma]) \text{ satisfy}\\
q_1, q_2 \notin \{ t_1,\ldots,t_n \}
\end{array}
 \right\}.
\end{equation*}
We have the extended map $ \widehat{M}^{0}_{{\it H}} \cup \widehat{M}^{10}_{{\it H}} \ra \mathrm{Hilb}^2(\cK'_5)$ by
\begin{equation*}
\widehat{M}^{10}_{{\it H}} \ni
(E,\Phi,\varphi,[\sigma]) \longmapsto (\{ (q,p) , (q,-p) \}, \lambda_+= \infty ) \in\mathrm{Hilb}^2(\cK'_5)
\end{equation*}
where $(q,p),(q,-p)$ are the pairs of the apparent singularities and the dual parameters, 
and $\lambda_+$ is the parameter of the fiber of the Hilbert-Chow morphism, that is, $p_2-p_1=\lambda_+(q_2-q_1)$.

Next, we extend the map $ \widehat{M}^{0}_{{\it H}} \cup \widehat{M}^{10}_{{\it H}} \ra \mathrm{Hilb}^2(\cK'_5)$ to $\widehat{M}_{{\it H}} \ra \mathrm{Hilb}^2(\cK'_5)$.
For $(E,\Phi, \varphi,[\sigma]) \in  \widehat{M}_{{\it H}}^1$, we describe the Higgs field $\Phi$ as follows:
\begin{equation*}
\Phi = 
\begin{cases}
 A_z^1 \otimes \omega_z    & \text{on } U_0 \\
 R_0^{-1}(A_z^1 \otimes \omega_z )R_0  & \text{on } U_{\infty}
\end{cases}\quad \text{where }\ 
A_z^1 :=\left(
\begin{array}{clcl}
f_{11}^{(3)}(z) & f_{12}^{(6)}(z) \\
 f_{21}^{(0)}(z) & -f_{11}^{(3)}(z)
\end{array}
\right).
\end{equation*}
By the automorphism of $E\cong \cO(1)\oplus \cO(-2)$, we can normalize $A_z^1$ as follows:
\begin{equation}\label{normalization O1+O-2 n=5}
A_z^1 =
\left(
\begin{array}{clclc}
0 & \tilde{f}_{12}^{(6)}(z) \\
 1 & 0
\end{array}
\right)
\end{equation}
where we put $\tilde{f}_{12}^{(6)}(z):=b_0+b_1z+\cdots+b_6z^6$.
The spectral curve $C_s \subset \cK_5$ is defined by $\eta^2 - f_{12}^{(6)}(z) =0$. 
The curve $C_s$ passes through the points $(t_1,\hat{\nu}_1), (t_1,-\hat{\nu}_1) , \ldots, (t_n,\hat{\nu}_n), (t_n,-\hat{\nu}_n)$, and $(q,p),(q,-p)$.
Here, $(q,p),(q,-p)$ are the pairs of the apparent singularities and the dual parameters.
We define the map $\widehat{M}^{1}_{{\it H}} \setminus \widehat{M}^{10}_{{\it H}}\ra \mathrm{Hilb}^2(\cK'_5)$ by
\begin{equation*}
\widehat{M}^{1}_{{\it H}}\ni (E,\Phi,\varphi,[\sigma]) \setminus \widehat{M}^{10}_{{\it H}} 
\longmapsto \{ (t_i,\hat{\nu}_i, v ), (t_i,-\hat{\nu}_i,-v) \}  \in \mathrm{Hilb}^2(\cK'_5)
\end{equation*}
where 
$\{ (t_i, \hat{\nu}_i), (t_i, -\hat{\nu}_i) \}$ are the apparent singularities and the dual parameters of $(E,\Phi,\varphi,[\sigma])$, and 
\begin{equation*}
v= \lim_{q\ra t_i} \frac{p  - \hat{\nu}_i }{q - t_i}= \frac{1}{ 2 \hat{\nu}_i} \frac{d}{dz}\tilde{f}_{12}^{(6)} (t_i).
\end{equation*}
Then we have the natural extended map $ \widehat{M}_{{\it H}} \ra \mathrm{Hilb}^2(\cK'_5)$.

Let $\widetilde{\mathrm{Hilb}}^2(\cK'_5)$ be the blowing-up of $\mathrm{Hilb}^2(\cK'_5)$ along $Z$.
We show that the spectral curves are determined by the point of $\widetilde{\mathrm{Hilb}}^2(\cK'_5)$.
Let $\tilde{x}:=( \{ (q_1,p_1, v_1^{i,\pm}), (q_2,p_2, v_2^{i,\pm}) \}, \lambda_+, \lambda_- )$ 
be a point of $\widetilde{\mathrm{Hilb}}^2(\cK'_5)$ where $p_2-p_1 = \lambda_+(q_2-q_1)$ and $p_2+p_1 = \lambda_-(q_2-q_1)$.
First, any spectral curves pass through the points $(t_i, \hat{\nu}_i)$ and $(t_i, -\hat{\nu}_i)$ for $i=1,\ldots,n$, that is, the polynomial $f_{12}^{(6)}(z)$ satisfies the condition
$\hat{\nu}_i^2  - f_{12}^{(6)}(t_i)=0$.
By the equations, we can determine the coefficients $b_1,b_2,b_3,b_6$ by $b_4,b_5$.
Second, the spectral curves passes through the points $\{ (q_1,p_1), (q_2,p_2) \}$, that is, $p_i^2- f_{12}^{(6)}(q_i)=0$ for $i=1,2$.
If $q_1 \neq q_2$, then the coefficients $b_4$ and $b_5$ are determined by $(q_1,p_1)$ and $(q_2,p_2)$,
that is, the spectral curve is determined by the apparent singularities and the dual parameters.
We consider the behavior of the spectral curve as $q_2 \ra q_1$ and $p_2\ra -p_1$.
We consider the following equations
\begin{equation*}
\begin{cases}
p_1^2- \tilde{f}_{12}^{(6)}(q_1)=0\\
\left( -p_1+ \lambda_- (q_2-q_1) \right)^2 - \tilde{f}_{12}^{(6)}(q_2)=0.
\end{cases}
\end{equation*}
as $q_2 \ra q_1$ and $p_2\ra -p_1$.
When $q_1 , q_2 \notin \{ t_1,\ldots,t_5\}$, we can determine the coefficients $b_4$ and $b_5$ by $q_1,p_1$ and $\lambda_-$.
We consider the case $q_1 =q_2 = t_i$ for some $i$.
Let $\lambda_-^i$ be the parameter such that $v_2^{i,+} + v_1^{i,+}=\lambda^i_- (q_2 -q_1)$, which is a blowing-up parameter of $\widetilde{\mathrm{Hilb}}^2(\cK'_5) \ra \mathrm{Hilb}^2(\cK'_5)$.
When we take $q_1 \ra t_i$ and $q_2 \ra t_i$,
we have the following equations
\begin{equation*}
\begin{cases}
 v_1^{i,+} = \frac{1}{ 2 \hat{\nu}_i} \frac{d}{dz} \tilde{f}_{12}^{(6)} (t_i)\\
\lambda_- = \frac{1}{2} \frac{1}{ 2 \hat{\nu}_i} \left( \frac{d^2}{dz^2}\tilde{f}_{12}^{(6)} (t_i)- \frac{1}{2\hat{\nu}_i^2} (\frac{d}{dz}\tilde{f}_{12}^{(6)} (t_i))^2 \right).
\end{cases}
\end{equation*}
By these equations, we can determine the coefficients $b_4$ and $b_5$ by $q_1,v_1^{i,+}$ and $\lambda_-$.
Therefore, spectral curves are determined by the point of $\widetilde{\mathrm{Hilb}}^2(\cK'_5)$.
If the underlying vector bundles of Higgs bundles are $\cO(1)\oplus \cO(-2)$,
then Higgs fields are determined by points of $\widetilde{\mathrm{Hilb}}^2(\cK'_5)$ by the normalization (\ref{normalization O1+O-2 n=5}).
On the other hand, cyclic vectors of bundle type $\cO(1)\oplus \cO(-2)$ are determined by the apparent singularities $q_1(=q_2)$. 
Then we have that the restriction map $ \widehat{M}^1_{{\it H}} \ra \mathrm{Hilb}^2(\cK'_5)$ is injective.
Finally, we obtain that the restriction map $ \widehat{M}_{{\it H}} \ra \mathrm{Hilb}^2(\cK'_5)$ is injective.
\end{proof}

We can describe the image of the map (\ref{equ of main theorem +}) as follows.
We define parameters $(v^{j,\pm})_{1\le j\le 5}$ by
\begin{equation*}
p -\nu_j^{\pm}(t_j-t_{k_1})(t_j-t_{k_2})(t_j-t_{k_3}) = v^{j,\pm} (q-t_j).
\end{equation*}
The parameters $(v^{j,\pm})_{1\le j\le 5}$ are blowing-up parameters of $\mathrm{Bl}_{\hat{\nu}^{\pm}_i} \mathbb{L} \ra \mathbb{L}$.
The moduli space $\widehat{M}_{{\it H}}$ is stratified as $\widehat{M}_{{\it H}} = \widehat{M}^{0}_{{\it H}} \cup \widehat{M}^{1}_{{\it H}}$
where $\widehat{M}^{i}_{{\it H}}$ is the locus such that $(E,\Phi,\varphi,[\sigma]) \in \widehat{M}^{i}_{{\it H}}$ satisfies $E \cong \cO(i)\oplus\cO(-i-1)$.
Then the images of $\widehat{M}^{0}_{{\it H}}$ and $\widehat{M}^{1}_{{\it H}}$ in $\widetilde{\mathrm{Hilb}}^2(\cK'_5)$ are the following
\begin{equation}\label{explicit description M0}
\begin{aligned}
\widehat{M}^{0}_{{\it H}} \cong&\  \{ ( \{ (q_1, p_1, v_1^{i,\pm}), (q_2,p_2, v_2^{i,\pm}) \}  ) \mid  q_1\neq q_2,v_j^{i,\pm}\in\mathbb{C}  \}\cup \\
&  \bigcup \left\{ ( \{ (q_1, p_1, v_1^{i,\pm}), (q_2,p_2, v_2^{i,\pm}) \} ,\lambda_+ )\ \middle|\  
\begin{array}{l}
q_1, q_2\notin \{ t_1,\ldots,t_5\}  ,\\
q_1 - q_2=p_1-p_2=0, \lambda_+\in \mathbb{C}
\end{array}
\right\}\cup \\
&  \bigcup_{\stackrel{j=1,\ldots,5}{\epsilon=\pm }}
 \left\{ ( \{ (q_1, p_1, v_1^{i,\pm}), (q_2,p_2, v_2^{i,\pm}) \},\lambda_+ , \lambda^{i}_+ )\ \middle| \  
\begin{array}{l}
 q_1 = q_2=t_j,\\
p_1=p_2=\nu_j^{\epsilon}(t_j-t_{k_1})(t_j-t_{k_2})(t_j-t_{k_3})\\
 v_1^{j,\epsilon}=v_2^{j,\epsilon} =\lambda_+\in \mathbb{C}
\end{array}
\right\}\\
\end{aligned}
\end{equation}
\begin{equation}\label{explicit description M1}
\begin{aligned}
\widehat{M}^{1}_{{\it H}} \cong &
\left\{ ( \{ (q_1, p_1, v_1^{i,\pm}), (q_2,p_2, v_2^{i,\pm}) \} ,\lambda_+, \lambda_- ) \ 
\middle|\  
\begin{array}{l}
q_1, q_2\neq t_j \text{ for any $j=1,\ldots,5$} ,\\
q_1-q_2= p_1+p_2=0,\\
 \lambda_+= \infty,\lambda_-\in \mathbb{C} 
\end{array}
 \right\}\cup\\
&\bigcup_{\stackrel{j=1,\ldots,5}{\epsilon=\pm }} \left\{ ( \{ (q_1, p_1, v_1^{i,\pm}), (q_2,p_2, v_2^{i,\pm}) \} ,\lambda_+,\lambda_-, \lambda_-^i ) \ 
\middle|\  
\begin{array}{l}
 q_1 = q_2=t_j,\\
p_1=-p_2=\nu_j^{\epsilon}(t_j-t_{k_1})(t_j-t_{k_2})(t_j-t_{k_3})\\
\lambda_+=\infty\\
 v_1^{j,\epsilon}=-v_2^{j,-\epsilon} =\lambda_-\in \mathbb{C}
\end{array}
 \right\}.
 \end{aligned}
\end{equation}
Here, $\lambda_{\pm}$ and $\lambda^j_{\pm}$ satisfy the following relations
\begin{equation*}
p_1-p_2=\lambda_+ (q_1-q_2),\quad  v_1^{i,\epsilon} - v_2^{i,\epsilon}= \lambda^{i}_+ (q_1-q_2 )  
\end{equation*}
\begin{equation*}
p_1+p_2 = \lambda_-(q_1-q_2),\quad v_1^{i,\epsilon} + v_2^{i,\epsilon}= \lambda^{i}_- (q_1-q_2 )  .
\end{equation*}
\begin{Rem}
We consider the map
\begin{equation*}
\begin{aligned}
\widehat{M}^{1}_{{\it H}} &\lra \mathbb{P}^1\\
(E,\Phi,\varphi,[\sigma]) &\longmapsto q_1
\end{aligned}
\end{equation*}
given by the description (\ref{explicit description M1}) and the natural projection.
Any fiber of this map is $\mathbb{C}^2$, which is isomorphic to $M_{{\it H}}^1$.
Note that $q_1$ is a coordinate of $\mathbb{P}H^0(\mathbb{P},E)$, which is a blowing-up parameter of $\widehat{M}_{{\it H}} \ra M_{{\it H}}$.
\end{Rem}

\section{Jumping families for Higgs bundles}\label{universel family}

In this section, we give an explicit description of the universal family of the moduli space $\widehat{M}_{{\it H}}$ for $n=5$.
For the purpose, we need give a description of jumping family,
which is a family of Higgs fields such that for generic parameters, the underlying vector bundles are $\cO\oplus \cO(-1)$ and
for special parameters, the underlying vector bundles are $\cO(1)\oplus \cO(-2)$.
Descriptions of jumping families are given by the lower and upper modifications.
In \ref{family for n}, we apply the description of jumping families to the case $n\ge 4$.
Then we obtain explicit descriptions of jumping families for the case $n\ge 4$.

\subsection{Jumping family for $n=5$}\label{construction of the family on cV1}

Suppose that $n=5$.
We consider the following covering of $\widehat{M}_{{\it H}}$:
\begin{equation*}
\begin{aligned}
\cV_0:=&\  \widehat{M}^{0}_{{\it H}} \subset \widehat{M}_{{\it H}} ,\\
\cV_1:= &\ 
\left\{ (E,\Phi,\varphi) \in M^{00}_{{\it H}}\ \middle|
\begin{array}{l}
\{ \tilde{p}_1 , \tilde{p}_{2} \} \not\subset \mathrm{Sym}^2 [\text{ branch points of $C_s$ }]  
\end{array}
  \right\}
\cup \widehat{M}^{1}_{{\it H}} \subset \widehat{M}_{{\it H}},
\end{aligned}
\end{equation*}
where $C_s$ is the spectral curve of $(E,\Phi,\varphi)$
and $\tilde{p}_{i}=(q_i,p_i)\in C_s$ is an apparent singularity and its dual of $(E,\Phi,\varphi)$.
Since we consider $\mathfrak{sl}_2$-Higgs bundles, $p_i=0$ implies that $\tilde{p}_{i}=(q_i,p_i)\in C_s$ is a branch point of $C_s$.
By Theorem \ref{Prop M^0 to Hilb}, we have an explicit description of the universal family $(E_{\cV_0},\Phi_{\cV_0})$ on $\cV_0\times \mathbb{P}^1$.
Now we give an explicit description of the universal family $(E_{\cV_1},\Phi_{\cV_1},[\sigma_{\cV_1}])$ on $\cV_1\times \mathbb{P}^1$.
Then we have an explicit description of the universal family $(\tilde{E},\tilde{\Phi},[\tilde{\sigma}])$ on $ \widehat{M}_{{\it H}}  \times \mathbb{P}^1$:
\begin{equation*}
\xymatrix{
(E_{\cV_1},\Phi_{\cV_1}, [\sigma_{\cV_1}]) \ar[r]^-{\subset} \ar[d] & (\tilde{E},\tilde{\Phi},[\tilde{\sigma}]) \ar[d] &  (E_{\cV_0},\Phi_{\cV_0}) \ar[d] \ar[l]_-{\supset} \\ 
\cV_1\times \mathbb{P}^1 \ar[r]^-{\subset} & \widehat{M}_{{\it H}} \times \mathbb{P}^1 &  \cV_0\times \mathbb{P}^1 \ar[l]_-{\supset} \rlap{.}
}
\end{equation*}
Let
\begin{equation}\label{family diag}
\Phi_{\cV_0} = 
\begin{cases}
 A_z^0 \otimes \omega_z    & \text{on } U_0 \\
 R_0^{-1}(A_z^0 \otimes \omega_z )R_0  & \text{on } U_{\infty}
\end{cases},\ 
\text{ where } 
A_z^0= \left(
\begin{array}{clcl}
a_1 z+a_0 & f_{12}^{(4)}(z) \\
 (z-q_1)(z-q_2) & -(a_1z+a_0)
\end{array}
\right),
\end{equation}
be the family on $\cV_0$ obtained by Theorem \ref{Prop M^0 to Hilb} for $n=5$.
Here we set
\begin{equation*}
a_1 := \frac{p_1-p_2}{q_1-q_2}, \quad a_0:=-\frac{p_1q_2 - p_2q_1}{q_1-q_2}, \quad 
R_0=
\begin{pmatrix}
1 & 0 \\
0 & 1/z 
\end{pmatrix}
\end{equation*}
and
we assume that $q_1,q_2\neq \infty$ for simplicity.
Set 
\begin{equation}\label{X parame space}
X:= \left\{ ((q_1,p_1), (q_2,p_2) , \lambda) \in (\cK_5')^2 \times \mathbb{C}\ \middle|
\begin{array}{l}
 p_2-p_1 = \lambda (q_2-q_1),  \\
p_1\neq 0, \text{and } q_2-q_1\neq0
\end{array}
  \right\}
\end{equation}
and 
\begin{equation}\label{hatX parame space}
\widehat{X}:= X \cup \left\{ ((q_1, p_1), (q_2, p_2) , \lambda ) \in (\cK_5')^2 \times \mathbb{C}\ \middle|
\ p_2-p_1 = q_2-q_1=0  \right\}.
\end{equation}

Let $P_1$, $P_2$, and $P_3$ be the following matrices
\begin{equation}\label{matrces P1P2P3} 
P_1:=\left(
\begin{array}{ll}
 1 & 0 \\
 0 & z-q_1
\end{array}
\right),\ 
P_2:=\left(
\begin{array}{ll}
 1 & 0 \\
 \frac{q_1-q_2}{2 p_1} & 1
\end{array}
\right),\ 
P_3:=\left(
\begin{array}{ll}
 \frac{1}{z-q_1} & 0 \\
 0 & 1
\end{array}
\right).
\end{equation}

\begin{Prop}\label{Jump descrip Prop}
\textit{
We define a family of $\bnu$-$\mathfrak{sl}_2$-parabolic Higgs bundle $(E_X, \Phi_X, \varphi_X, [\sigma_X])$ with a cyclic vector
on $X \times \mathbb{P}^1$ as 
\begin{equation}\label{family diag jump}
\begin{aligned}
\Phi_{X} =
\begin{cases}
 (P_1 P_2 P_3)^{-1}(A^0_z\otimes \omega_z) P_1 P_2 P_3     & \text{on } U_0 \\
 R_0^{-1}(A_z^0 \otimes \omega_w )R_0  & \text{on } U_{\infty}^{q_1}
\end{cases}
\end{aligned}
\end{equation}
where $U_{\infty}^{q_1} = \Spec\bC\left[w,1/(q_1w-1)\right]$.
Here $\sigma_X$ is the element of $H^0( \mathbb{P}^1, E_X)$ such that the zero of $\sigma_X$ is $q_1$ when $E_X \cong \cO(1) \oplus \cO(-2)$.
We can extend this family on $X \times \mathbb{P}^1$ to the family on $\widehat{X} \times \mathbb{P}^1$, naturally.
For the extended family, we have the following.
If $q_1\neq q_2$, then the underlying vector bundle is $\cO \oplus \cO(-1)$.
If $q_1=q_2$, then the underlying vector bundle is $\cO(1) \oplus \cO(-2)$.
}
\end{Prop}

\begin{proof}
We describe the construction of a family of $\bnu$-$\mathfrak{sl}_2$-parabolic connections with a cyclic vector
on $X \times \mathbb{P}^1$ so that this family satisfies the assertions of the proposition.
As the result, we obtain the family which has the description (\ref{family diag jump}).
Then this proposition follows from this construction.

By the natural map $X \ra \cV_0$, the family $\Phi_{\cV_0}$ on $\cV_0 \times \mathbb{P}^1$ induces the family on $X \times \mathbb{P}^1$, denoted by $\Phi_X^0$.
We consider the lower and upper modifications
$(q_1,l_{p_1})^{\text{up}} \circ (q_1,l_{p_1})^{\text{low}}  (\Phi_X^0) $, denoted by $\Phi_X$,
where $l_{p_1}$ is a one dimensional subspace of $E|_{q_1}$ which corresponds to the eigenspace of the residue of $\Phi_X^0$ at $q_1$ with the eigenvalue $p_1$.
Explicitly, the modifications are described as follows.
We consider the following diagram
\begin{equation*}
\xymatrix{
U_0 \times \bC^2  & U_{\infty}^{q_1}\times \bC^2 \ar[l]_-{R_0} \ar[ddl]^-{T(q_2)} \\
U_0 \times \bC^2   \ar[u]^-{P_1} \\ 
U_0 \times \bC^2   \ar[u]^-{P_2\circ P_3}
}
\end{equation*}
where 
\begin{equation*}
T:=\left(
\begin{array}{ll}
 z-q_1 & 0 \\
  -\frac{q_1-q_2}{2 p_1} & \frac{1}{z(z-q_1)}
\end{array}
\right).
\end{equation*}
In particular,
\begin{equation}\label{trans O1+O-2}
T^0_{q_1}=\lim_{q_2\rightarrow q_1}T=
\left(
\begin{array}{ll}
 z-q_1 & 0 \\
 0 & \frac{1}{z(z-q_1)}
\end{array}
\right).
\end{equation}
Here, the transformation $P_1$ implies the lower modification $(q_1,l_{p_1})^{\text{low}}$
and the transformation $P_2 \circ P_3$ implies the upper modification $(q_1,l_{p_1})^{\text{up}}$.
Namely, we describe $\Phi_{X}$ as
\begin{equation}
\begin{aligned}
\Phi_{X} = &\,\,(q_1,l_{p_1})^{\text{up}} \circ (q_1,l_{p_1})^{\text{low}}  (\Phi_{X})\\
= 
&\begin{cases}
 (P_1 P_2 P_3)^{-1}(A^0_z\otimes \omega_z) P_1 P_2 P_3      & \text{on } U_0 \\
 R_0^{-1}(A_z^0 \otimes \omega_w )R_0  & \text{on } U_{\infty}^{q_1}.
\end{cases}
\end{aligned}
\end{equation}
Taking the limit $q_2\ra q_1$, we have Higgs bundles of bundle type $\cO(1)\oplus\cO(-2)$.
Then we have the description of the family of the Higgs fields (\ref{family diag jump}), 
and the family satisfies the assertion of the proposition.
By the transition function (\ref{trans O1+O-2}), the zero of cyclic vectors is $q_1$ when $E_X \cong \cO(1) \oplus \cO(-2)$.
\end{proof}
By this proposition, we have a map $\widehat{X} \ra \cV_1\subset \widehat{M}_H$.

\begin{figure}[htbp]
 \begin{center}
  \includegraphics[width=80mm]{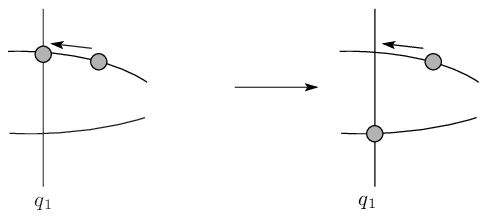}
 \end{center}
 \caption{The apparent singularities and the dual parameters}
 \label{fig2}
\end{figure}

Next we compute the apparent singularities and the dual parameters of $\Phi_{X}$ when $q_1\neq q_2$.
Put 
\begin{equation}\label{matrices Q1Q2}
Q_1:=
\left(
\begin{array}{ll}
 z-q_1 & \frac{2 p_1}{q_1-q_2} \\
 -\frac{q_1-q_2}{2 p_1} & 0
\end{array}
\right), \quad
Q_2:=
\left(
\begin{array}{ll}
 1 & -\frac{2 p_1 w^2}{(q_1-q_2) (wq_1-1) } \\
 0 & 1
\end{array}
\right).
\end{equation}
Then we have the following diagram
\begin{equation*}
\xymatrix@=17pt{
U_0 \times \bC^2    &  U_{\infty}^{q_1}\times \bC^2  \ar[l]_-{T(q_2)} \\
U_0 \times \bC^2  \ar[u]^-{Q_1}   &  U_{\infty}^{q_1}\times \bC^2  \ar[l]^-{R_0}  \ar[u]_-{Q_2}\rlap{.}
}
\end{equation*}
By the transformation by $Q_1$ and $Q_2$,
we can describe $\Phi_{X}$ as in the description (\ref{normal form of O+O(-1)}). 
The apparent singularities and the dual parameters of $\Phi_{X}$ is $\{ (q_1 ,-p_1),(q_2, p_2)\}$ (Figure \ref{fig2}).
Then we have a map $X \ra \widetilde{\mathrm{Hilb}}^2(\cK'_5)$.
Moreover, we can extend this map to $\widehat{X} \ra \widetilde{\mathrm{Hilb}}^2(\cK'_5)$ by
\begin{equation*}
\widehat{X} \setminus X \ni ((q_1,p_1),(q_1,p_1),\lambda) \longmapsto ((q_1,-p_1),(q_1,p_1), \infty ,\lambda) \in \widetilde{\mathrm{Hilb}}^2(\cK'_5)
\end{equation*}
where $(\infty, \lambda)$ are values of a parameter of pre-images of the Hilbert-Chow morphism and of a blowing-up parameter of 
$\widetilde{\mathrm{Hilb}}^2(\cK'_5)\ra {\mathrm{Hilb}}^2(\cK'_5)$.
We obtain the commutative diagram
\begin{equation*}
\xymatrix@=15pt{
  \widehat{X} \ar[dr] \ar[d]   & \\
\cV_1 \ar@{^{(}->}[r]  & \widetilde{\mathrm{Hilb}}^2(\cK'_5).
}
\end{equation*}
By this diagram, 
the explicit description of the family (\ref{family diag jump}) induces
an explicit description of the universal family on $\cV_1 \times \mathbb{P}^1$ which is parametrized by 
the apparent singularities and their duals.

\subsection{Jumping families for $n\ge 5$}\label{family for n}
In this section, we give an explicit description of  jumping families for $n\ge5$ as in the previous section.

Let
\begin{equation}\label{family diag 0}
\Phi_{0} = 
\begin{cases}
 A_z^0 \otimes \omega_z    & \text{on } U_0 \\
 R_0^{-1}(A_z^0 \otimes \omega_w )R_0  & \text{on } U_{\infty}
\end{cases},\ 
\text{ where } 
A^0_z= \left(
\begin{array}{clcl}
f_{11}^{(n-2)}(z) & f_{12}^{(n-1)}(z) \\
 (z-q_1)\cdots (z-q_{n-3}) & -f^{(n-2)}_{11}(z)
\end{array}
\right),
\end{equation}
be the family on $M^0_{{\it H}}$ obtained by Theorem \ref{Prop M^0 to Hilb}.
We assume that $q_1,\ldots,q_{n-3}\neq \infty$ for simplicity.

First, we construct a family having Higgs bundles of bundle type $\cO(1)\oplus \cO(-2)$ from the family $\Phi_{0}$ by lower and upper modifications as in \ref{construction of the family on cV1}.
Fix $(q_1,p_1) \in \cK'_n$ and assume that $q_1 \neq q_2$. 
Put
\begin{equation*}
\Phi_{1}:= (q_1,l_{p_1})^{\text{up}} \circ (q_1,l_{p_1})^{\text{low}}  (\Phi_{0}).
\end{equation*} 
If we take the limit $q_2\ra q_1$ of $\Phi_{1}$, then we have Higgs bundles of bundle type $\cO(1)\oplus\cO(-2)$:
\begin{equation}\label{family diag 1}
\lim_{q_2\ra q_1} \Phi_{1} = 
\begin{cases}
 A_z^1 \otimes \omega_z    & \text{on } U_0 \\
 (R_1^{q_1})^{-1}(A_z^1 \otimes \omega_w )R_1^{q_1}  & \text{on } U^{q_1}_{\infty},
\end{cases} \quad \text{ where }
A_z^1= \left(
\begin{array}{clcl}
{}'f^{(n-2)}_{11}(z) & {}'f_{12}^{(n+1)}(z) \\
 (z-q_3)\cdots (z-q_{n-3}) & -{}'f^{(n-2)}_{11}(z)
\end{array}
\right)
\end{equation}
and 
\begin{equation*}
R^{q_1}_1= \left(
\begin{array}{clcl}
z-q_1 & 0 \\
 0 & \frac{1}{z(z-q_1)}
\end{array}
\right).
\end{equation*}
Here, the entries of $A_z^1$ satisfy the following equations:
\begin{equation*}
\begin{cases}
p_i- {}'f^{(n-2)}_{11}(q_i)=0 \ \text{ for $i=3,\ldots,{n-3}$,}\\
p_1^2- \det (A_z^1|_{z=q_1}) =0 ,& \\
\det \left( \res_{z=t_i} \lim_{q_2\ra q_1} \Phi_{1} \right) -\nu_i^2  =0\ \text{ for $i=1,\ldots,n$.}
\end{cases}
\end{equation*}
Note that $q_1$ is the zero of the corresponding cyclic vector $\sigma \in H^0(\bP^1, E)$.

Next, we construct a family having Higgs bundles of bundle type $\cO(2)\oplus \cO(-3)$ from the family $\Phi_{1}$.
Here, we assume that $q_1,\ldots,q_{n-3}\neq \infty$ for simplicity.
For the Higgs field (\ref{family diag}),
fix $(q_3,p_3) \in \cK'_n$ and assume that $q_3\neq q_4$.
Put
\begin{equation*}
\Phi_{2}:= (q_3,l_{p_3})^{\text{up}} \circ (q_3,l_{p_3})^{\text{low}}  \left(\lim_{q_2\ra q_1} \Phi_{1} \right)
\end{equation*} 
If we take the limit $q_4\ra q_3$ of $\Phi_{2}$, then we have Higgs bundles of bundle type $\cO(2)\oplus\cO(-3)$:
\begin{equation}\label{family diag 2}
\lim_{q_2\ra q_1} \Phi_{2} = 
\begin{cases}
 A_z^2 \otimes \omega_z    & \text{on } U_0 \\
 (R_2^{q_1})^{-1}(A_z^2 \otimes \omega_w )R_2^{q_1}  & \text{on } U^{q_1}_{\infty}
\end{cases} \quad \text{ where }
A_z^2= \left(
\begin{array}{clcl}
{}''f^{(n-2)}_{11}(z) & {}''f_{12}^{(n+3)}(z) \\
 (z-q_5)\cdots (z-q_{n-3}) & -{}''f^{(n-2)}_{11}(z)
\end{array}
\right)
\end{equation}
and 
\begin{equation*}
R^{q_1}_2= \left(
\begin{array}{clcl}
(z-q_1)(z-q_3) & 0 \\
 0 & \frac{1}{z(z-q_1)(z-q_3)}
\end{array}
\right).
\end{equation*}
Here, the entries of $A_z^2$ satisfy the following equations:
\begin{equation*}
\begin{cases}
p_i- {}''f^{(n-2)}_{11}(q_i)=0\ \text{ for $i=5,\ldots,{n-3}$,}\\
p_i^2- \det (A_z^2|_{z=q_i}) =0\ \text{ for $i=1,3$,} \\
\det \left(\res_{z=t_i} \lim_{q_2\ra q_1} \Phi_{2} \right) -\nu_i^2 =0\ \text{ for $i=1,\ldots,n$.}\\
\end{cases}
\end{equation*} 
Note that $q_1,q_3$ is the zeros of the corresponding cyclic vector $\sigma \in H^0(\bP^1, E)$.
We continue this process.
Then we have family having Higgs bundles of bundle type $\cO(k)\oplus \cO(-k-1)$ for $k=1,\ldots,[(n-3)/2]$.

\section{Geometric description for connection cases}\label{connection case}

Suppose that $\bnu$ satisfies the condition (\ref{generic condition 2}) and $\nu_1\cdots\nu_n\neq0$.
We put 
\begin{equation*}
\begin{aligned}
(t_1,\ldots,t_{n})&:=(0,1,x_1,\ldots,x_{n-3},\infty),\\
(\nu^{\pm}_1,\ldots,\nu_{n-1}^{\pm},\nu^{+}_n,\nu_n^{-})&:=(\pm\nu_0,\pm\nu_1,\ldots,\pm\nu_{n-1},\nu_{n},1-\nu_{n}), \text{ and}\\
\hat{\nu}_i &:=  \nu_i (t_i-t_1)\cdots(t_i-t_{i-1})(t_i-t_{i+1})\cdots (t_i-t_{n-1}) \text{ for $i=1,\ldots,n-1$}.
\end{aligned}
\end{equation*}
Let $M^k$ (resp. $\widehat{M}^k$) be the subvariety of $M$ (resp. $\widehat{M}$) where $E\cong \cO(k) \oplus\cO(-k-1)$.
First, we define the apparent singularities and the dual parameters of $(E,\nabla,\varphi) \in M^0$, and
we compute the apparent singularities and the dual parameters for $n\ge4$.
Then we have Theorem \ref{from M0con to Hilb 0} (Theorem \ref{from M0con to Hilb}). 
Second, we assume $n=5$.
We construct a jumping family.
Third, we compute the apparent singularities of this jumping family on the locus of bundle type $\cO\oplus \cO(-1)$,
and we analyze the behavior of the apparent singularities of the jumping families when the parameter closes to the jumping locus.
Then we obtain a map from $\widehat{M}$ to the Hilbert scheme of points on some surface.
Moreover, we take some sequence of blowing-ups of the Hilbert scheme.
Then we obtain an injective map from $\widehat{M}$ to the blowing-ups.

\subsection{Geometric description of $M^0$ for $n\ge 4$}\label{connection O+O(-1)}

Let $(E,\nabla,\varphi) \in M$.
We can define the \textit{apparent singularities of $(E,\nabla,\varphi) \in M$} as follows.
We fix a section $\sigma \in H^0(\bP^1, E)$.
For the section $\sigma$, we define the following composition
\begin{equation*}
\cO_{\bP^1} \xrightarrow{\ \sigma\ } E \xrightarrow{\ \nabla\ } E \otimes L \lra (E/\cO_{\bP^1}) \otimes L.
\end{equation*}
The composition $\cO_{\bP^1}\ra (E/\cO_{\bP^1}) \otimes L$ is an $\cO_{\bP^1}$-morphism, which is injective.
Then we can define a subsheaf $F^0\subset E $ such that $\cO_{\bP^1} \ra (F^0/\cO_{\bP^1}) \otimes L$ is an isomorphism.
By the isomorphism $F^0/\cO_{\bP^1} \cong L^{-1}$, we have $F^0 \cong \cO_{\bP^1} \oplus  L^{-1}$.
Therefore, we have the following exact sequence
\begin{equation}\label{ES of App for conn}
0 \lra \cO_{\bP^1} \otimes L^{-1} \lra E \lra T_A \lra 0
\end{equation}
where $T_A$ is a torsion sheaf.
By the Riemann-Roch theorem, we have that the torsion sheaf $T_A$ is length $n-3$.

\begin{Def}
For $(E,\nabla , \varphi) \in M$ and a nonzero section $\sigma \in H^0(\bP^1,E)$, 
we call the support of $T_A$ {\it apparent singularities of a $\bnu$-$\mathfrak{sl}_2$-parabolic connection with a cyclic vector $(E,\nabla , \varphi,[\sigma])$}.
\end{Def}

For $(E,\nabla , \varphi) \in M^0$, we define \textit{dual parameters} as follows.
Since $E\cong \cO\oplus \cO(-1)$, we can denote the connection $\nabla$ by 
\begin{equation*}
\nabla = 
\begin{cases}
d+ A_z^0 \otimes \omega_z    & \text{on } U_0 \\
d+ R_0^{-1}dR_0 +  R_0^{-1}(A_z^0 \otimes \omega_z )R_0  & \text{on } U_{\infty}
\end{cases}\quad \text{where }\ 
A_z^0 :=\left(
\begin{array}{clcl}
f_{11}^{(n-2)}(z) & f_{12}^{(n-1)}(z) \\
 f_{21}^{(n-3)}(z) & -f_{11}^{(n-2)}(z)
\end{array}
\right).
\end{equation*}
Note that the zeros of the polynomial $f_{21}^{(n-3)}(z)$ are the apparent singularities of $(E,\nabla , \varphi)$.
We denote by $\{ q_1,\ldots,q_{n-3} \}$ the apparent singularities.
We put $p_i := f_{11}^{(n-2)}(q_i) \in L|_{q_i}$.
We call $\{ p_1 ,\ldots, p_{n-3} \}$ the \textit{dual parameters} of $(E,\nabla , \varphi) \in M^0$.
The definition of the apparent singularities and the dual parameters is already given by Oblezin \cite[Section 3]{Obl}.

Let $\widetilde{\cK}_n'$ be the Zariski open set of the blowing-up of Hirzebruch surface of degree $n-2$ defined in \ref{Hirzebruch blow up},
and $\widetilde{\cK}_n$ be the contraction $\widetilde{\cK}'_n\ra \widetilde{\cK}_n$. 
Then we can define the following map
\begin{equation}\label{M0 to Sym}
\begin{aligned}
 M_{}^0 &\lra \mathrm{Sym}^{n-3}(\widetilde{\cK}_n) \\
 (E,\nabla, \varphi) &\longmapsto \{(q_1,p_1),\ldots,(q_{n-3},p_{n-3}) \},
\end{aligned}
\end{equation}
which is already constructed in \cite[Section 3]{Obl}.
We consider the composite of the Hilbert-Chow morphism and the blowing-up
\begin{equation*}
\Hilb^{n-3}(\widetilde{\cK}'_n) \lra \mathrm{Sym}^{n-3}(\widetilde{\cK}'_n) \lra \mathrm{Sym}^{n-3}(\widetilde{\cK}_n).
\end{equation*}
where $\widetilde{\cK}_n' \ra \widetilde{\cK}_n$ is the blowing up defined in \ref{Hirzebruch blow up}.
By the same argument as in the proof of Theorem \ref{Prop M^0 to Hilb}, we obtain the following
\begin{Thm}\label{from M0con to Hilb}
\textit{
We can extend the map} (\ref{M0 to Sym}) \textit{to
\begin{equation*}
 M_{}^0 \lra \mathrm{Hilb}^{n-3}(\widetilde{\cK}'_n).
\end{equation*}
This map is injective.
Moreover, we can give an explicit description of the universal family $(\tilde{E}^{(0)}, \tilde{\nabla}^{(0)}) \ra  M^0\times \mathbb{P}^1$.
}
\end{Thm}

\subsection{Jumping family for $n=5$}\label{jumping phenomenon conn}
Suppose that $n=5$.
In this section, we give an explicit description of a jumping family of connections.

Let  $\cK_n'$ be the Zariski open set of the blowing-up of Hirzebruch surface of degree $n-2$ corresponding to the moduli space of parabolic Higgs bundles (defined in \ref{Hirzebruch blow up}),
and $\cK_n$ be the contraction $\cK'_n\ra \cK_n$. 
Fix $(e_0,e_1) \in \mathbb{C}^2$.
Set 
\begin{equation}
X:= \left\{ ((q_1,p_1), (q_2,p_2) , \lambda) \in (\cK_5')^2 \times \mathbb{C}\ \middle|
\begin{array}{l}
 p_2-p_1 = \lambda (q_2-q_1),  \\
p_1\neq 0, \text{and } q_2-q_1\neq0
\end{array}
  \right\}
\end{equation}
and 
\begin{equation}
\widehat{X}:= X \cup \left\{ ((q_1, p_1), (q_2, p_2) , \lambda ) \in (\cK_5')^2 \times \mathbb{C}\ \middle|
\ p_2-p_1 = q_2-q_1=0  \right\},
\end{equation}
which are defined in \ref{construction of the family on cV1}, and
let $Q_1, Q_2$ and $R_0$ be the following matrices
\begin{equation*}
Q_1:=
\left(
\begin{array}{ll}
 z-q_1 & \frac{2 p_1}{q_1-q_2} \\
 -\frac{q_1-q_2}{2 p_1} & 0
\end{array}
\right), \quad
Q_2:=
\left(
\begin{array}{ll}
 1 & -\frac{2 p_1 w^2}{(q_1-q_2) (wq_1-1) } \\
 0 & 1
\end{array}
\right), \text{ and }
R_0=
\begin{pmatrix}
1 & 0 \\
0 & 1/z 
\end{pmatrix},
\end{equation*}
respectively. ($Q_1$ and $Q_2$ was defined in \ref{construction of the family on cV1}). 
\begin{Prop}\label{Jump descrip Prop}
\textit{
We define a family of $\bnu$-$\mathfrak{sl}_2$-parabolic connections $(E_X,\nabla_X,\varphi_X, [\sigma_X])$ with a cyclic vector
on $X \times \mathbb{P}^1$ as 
\begin{equation}\label{jump family descrip conn}
\nabla_X=
\begin{cases}
d+ Q_1 d Q_1^{-1}+ Q_1 \widetilde{A}_{U_0} Q_1^{-1}    & \text{on } U_0 =\Spec\bC\left[z\right] \\
d+ Q_2 d Q_2^{-1} + Q_2 \left( R_0^{-1} d R_0+ R_0^{-1} \widetilde{A}_{U_0} R_0 \right) Q_2^{-1}  & \text{on } U_{\infty}^{q_1}=\Spec\bC\left[w,1/(q_1w-1)\right].
\end{cases}
\end{equation}
Here $\widetilde{A}_{U_0}$ is defined below, and 
$\sigma_X$ is the element of $H^0( \mathbb{P}^1, E_X)$ such that the zero of $\sigma_X$ is $q_1$ when $E_X \cong \cO(1) \oplus \cO(-2)$.
We can extend this family on $X \times \mathbb{P}^1$ to the family on $\widehat{X} \times \mathbb{P}^1$, naturally.
For the extended family, we have the following.
If $q_1\neq q_2$, then the underlying vector bundle is $\cO \oplus \cO(-1)$.
If $q_1=q_2$, then the underlying vector bundle is $\cO(1) \oplus \cO(-2)$.
}
\end{Prop}

We define a family of connection matrices $\widetilde{A}_{U_0}$ on $X \times U_0$ as 
\begin{equation}\label{connection matrix for jump}
\widetilde{A}_{U_0}=
\left(
\begin{array}{ll}
 F_{11}(z) & F_{12}(z) \\
 F_{21}(z) & -F_{11}(z)
\end{array}
\right) \frac{dz}{z(z-1)(z-x_1)(z-x_2)}
\end{equation}
where 
\begin{equation*}
\begin{aligned}
F_{11}(z):=&\ \frac{2p_1(z-q_1)}{q_2-q_1}+ (z-q_1)\lambda + \frac{1}{2} \biggl( (1 - 2 \nu_{5} e_0 (q_2-q_1)) z^3 +( - 2 \nu_{5} e_1 (q_2-q_1) + q_1 - x_1 - x_2 -1) z^2 \\
&  +\Bigl( -2(p_1 - \nu_{5} q_1^2 (q_2-q_1) ) e_0 + 2 \nu_{5} q_1 (q_2-q_1) e_1 -  \frac{\nu_{5}(q_2-q_1)}{p_1}  q_1(q_1-1) (q_1 - x_1) (q_1 - x_2) \\
&+ q_1^2  -  (1  + x_1 + x_2) q_1+ x_1 + x_2 + x_1 x_2\Bigr) z -2 p_1 (1 + e_1 + e_0 q_1) + (q_1-1) (q_1 - x_1) (q_1 - x_2) \biggr),\\ 
F_{21}(z):=&\ (1 + e_0 (q_1 - q_2))z^2 - (q_1 +q_2 + e_1( q_2 - q_1) ) z \\
&\qquad \qquad + q_1 q_2  + \frac{(q_1 - q_2)q_1}{2 p_1} (-2 e_1 p_1 - 2 e_0 p_1 q_1 + (q_1-1) (q_1 - x_1) (q_1 - x_2)).
\end{aligned}
\end{equation*}
We omit the description of $F_{12}(z)$, since the description is lengthened and is not necessary for computation of the apparent singularities and their duals.

By this proposition, we have a map $\widehat{X} \ra \widehat{M}$.
On the other hand, 
we have a map $X \rightarrow \mathrm{Hilb}^{n-3}(\widetilde{\cK}'_n)$ defined by
$ ( (q_1,p_1), (q_2,p_2), \lambda) \mapsto  (\{ (q_1',p_1'), (q_2',p_2')\}, \lambda)$
where $q_1'$, $q_2'$ are the zero of $F_{21}(z)$, and $p'_1, p_2'$ are their duals, that is, $p_i'=F_{11}(q_i')$.
Then we have the following diagram 
\begin{equation}\label{key diagram conn}
\xymatrix@=15pt{
  \widehat{X} \ar[d] \ar@{.>}[dr]   & \\
\widehat{M} \ar@{.>}[r]  & \mathrm{Hilb}^2(\widetilde{\cK}'_5).
}
\end{equation}

\begin{proof}[{Proof of Proposition \ref{Jump descrip Prop}}]
We describe the construction of a family of $\bnu$-$\mathfrak{sl}_2$-parabolic connections with a cyclic vector
on $X \times \mathbb{P}^1$ so that this family satisfies the assertions of the proposition.
As the result, we obtain the family which has the description (\ref{jump family descrip conn}).
Then this proposition follows from this construction.

Let $A_z^0 \otimes \omega_z$ be the Higgs field (\ref{family diag}) on $X \times U_0$ defined in \ref{construction of the family on cV1}.
Let $\Phi_{X}$ be the Higgs field on $X \times \mathbb{P}^1 $ defined by
\begin{equation}\label{representative of jump}
P^{-1}(A_z^0\otimes \omega_z) P    \text{ on } U_0, \quad \text{and} \quad  
 R_0^{-1}(A_z^0 \otimes \omega_z )R_0   \text{ on } U_{\infty}
\end{equation}
where the eigenvalues of the residue matrices are 
\begin{center}
  \begin{tabular}{|c|c|c|c|c|} \hline
    $\res_{0} \Phi_{X}$ & $\res_{1} \Phi_{X}$ & $\res_{t_1} \Phi_{X}$ & $\res_{t_2} \Phi_{X}$ & $\res_{\infty} \Phi_{X}$ \\ \hline 
   $\nu_1' ,- \nu_1'$  &  $\nu_2',- \nu_2' $ &$ \nu_3',- \nu_3' $&$ \nu_4',- \nu_4' $&$\nu_5',- \nu_5' $  \\ \hline
  \end{tabular}\ .
\end{center}
Here we put $P:=P_1 P_2 P_3$
where $P_1$, $P_2$, and $P_3$ are the matrices (\ref{matrces P1P2P3}).
Note that the underlying vector bundles have bundle type $\cO(1)\oplus \cO(-2)$ when $q_1= q_2$.
When $q_1 \neq q_2$, by $Q_1$ and $Q_2$ which are the matrices (\ref{matrices Q1Q2}), we denote the Higgs field $\Phi_{X}$ by
\begin{equation}\label{normal form of jumping family}
 Q_1^{-1}P^{-1}(A_z^0\otimes \omega_z) P Q_1     \text{ on } U_0, \quad \text{and} \quad
 Q_2^{-1} R_0^{-1}(A_z^0 \otimes \omega_z )R_0 Q_2  \text{ on } U_{\infty}
\end{equation}
as in the description (\ref{normal form of O+O(-1)}).

Now we construct a family of initial connections $\nabla_0$ on $X \times \mathbb{P}^1$ and we determine $\nu_1',\ldots,\nu_5'$
so that $\nabla_0 + \Phi_{X}$ is the desired family $\nabla_X$.
We put
\begin{equation}\label{T infty}
T_{\infty}=
\left(
\begin{array}{clclc}
1 & \nu_5' \\
 0 & 1
\end{array}
\right).
\end{equation}
Then we have
\begin{equation*}
T_{\infty}^{-1}  (\res_{t_n} \Phi_{\mathcal{V}_0}) T_{\infty} =
T_{\infty}^{-1}
\left(
\begin{array}{clclc}
0 & -\nu_5'^2 \\
 -1 & 0
\end{array}
\right)
T_{\infty}=
\left(
\begin{array}{clclc}
\nu_5' & 0 \\
-1 & -\nu_5'
\end{array}
\right).
\end{equation*}
Let
$\nabla_{0}\colon \cO\oplus \cO(-1) \ra (\cO\oplus \cO(-1)) \otimes \Omega^1_{\bP^1}(D)$
be the connection defined by 
\begin{equation}\label{nalba0 6.7.}
\begin{cases}
d+  R_{0,w}^{-1} d R_{0,w} + R_{0,w}^{-1} T_{\infty}(B_w \otimes \omega_w) T_{\infty}^{-1} R_{0,w}   & \text{on } U_0 \\
d+ T_{\infty}(B_w \otimes \omega_w) T_{\infty}^{-1}  & \text{on } U_{\infty}
\end{cases}
\end{equation}
where we put
$R_{0,w}:=\left(
\begin{array}{ll}
 1 & 0 \\
 0 & 1/w
\end{array}
\right)$ and  
\begin{equation*}
B_w:= \left(
\begin{array}{ll}
f_3w^3+f_2w^2 + f_1 w +f_0  & d_4w^4+d_3w^3+d_2w^2 +d_1w \\
(q_1-q_2)(e_2w^2+e_1w +e_0)  & (w-1)(x_1w-1)(x_2w-1) -(f_3w^3+f_2 w^2 + f_1 w +f_0)
\end{array}
\right)
\end{equation*}
where $f_0,\ldots,f_3,d_1,\ldots,d_4,e_0,e_1,$ and $e_2$ are parameters.
When $q_1\neq q_2$, we define the connection $\nabla_0 + \Phi_{X}$ as 
\begin{equation}\label{Def of jumping family}
\begin{cases}
d+ R_{0,w}^{-1} d R_{0,w} + R_{0,w}^{-1} T_{\infty}(B_w \otimes \omega_w) T_{\infty}^{-1} R_{0,w}+
 Q_1^{-1}P^{-1}(A_z^0\otimes \omega_z) P Q_1    & \text{on } U_0 \\
d+ T_{\infty}(B_w \otimes \omega_w) T_{\infty}^{-1}+ Q_2^{-1} R_0^{-1}(A_z^0 \otimes \omega_z )R_0 Q_2  & \text{on } U_{\infty}.
\end{cases}
\end{equation}
Then the eigenvalues of the residue matrix of $\nabla_0 + \Phi_{X}$ at $\infty$ is the following
\begin{center}
  \begin{tabular}{|c|c|} \hline
    $\res_{\infty} \Phi_{X}$ & $\res_{\infty} (\nabla_0 + \Phi_{X})$ \\ \hline 
    $\nu_5', -\nu_5'$ & $-f_0 + \nu_5', 1-(-f_0+\nu_5')$.  \\\hline
  \end{tabular}
\end{center}

We put $(\nabla_0)_{U_0}= R_{0,w}^{-1} d R_{0,w} + R_{0,w}^{-1} T_{\infty}(B_w \otimes \omega_w) T_{\infty}^{-1} R_{0,w}$.
First, we determine the parameters of $B_w$ so that the limit $\lim_{q_2\ra q_1}\left( Q_1 d(Q_1^{-1})+ Q_1 (\nabla_0)_{U_0} Q_1^{-1}\right)$ is convergence.
We claim that if we determine the parameters $f_0,\ldots,f_3$, $d_0,\ldots,d_4$, $e_1,e_2$, and $e_3$ such that 
the following polynomial 
\begin{equation*}
\begin{aligned}
&( 2 f_0 -1)z^4+ (1 + 2 f_1 - 2 f_0 q_1 + x_1 + x_2)z^3 +(2 f_2 + 2 e_0 p_1 - 2 f_1 q_1 - x_1 - x_2 - x_1 x_2)z^2\\
&\ \ +(2 f_3 + 2 e_1 p_1 - 2 f_2 q_1 + x_1 x_2)z+2 (e_2 p_1 - f_3 q_1)
\end{aligned}
\end{equation*}
is identically zero,
then the limit is convergence.
We solve the simultaneous linear equations as follows:
\begin{equation*}
\begin{aligned}
f_0&=\frac{1}{2}, \  f_1=\frac{1}{2}( q_1 - x_1 - x_2-1),\ f_2=\frac{1}{2} (-2 e_0 p_1 + q_1^2 + x_1 + x_2 + x_1 x_2 - q_1 (1 + x_1 + x_2)),\\
f_3&=\frac{1}{2} (-2 e_1 p_1 - 2 e_0 p_1 q_1 + (q_1 -1) (q_1 - x_1) (q_1 - x_2)),\text{ and}\\
e_2&=\frac{q_1}{2 p_1}  (-2 e_1 p_1 - 2 e_0 p_1 q_1 + (q_1 -1) (q_1 - x_1) (q_1 - x_2)).
\end{aligned}
\end{equation*}
Then we can define the limit $\lim_{q_2 \ra q_1} (\nabla_0+\Phi_{X})$ as 
\begin{equation*}
\begin{cases}
\lim_{q_2 \ra q_1}( d+ Q_1 d(Q_1^{-1})+ Q_1 (\widetilde{\nabla}_0)_{U_0} Q_1^{-1}+P(q_2)^{-1}(A_z^0\otimes \omega_z) P(q_2))    & \text{on } U_0 \\
\lim_{q_2 \ra q_1}(d+Q_2 d(Q_2^{-1})+ Q_2 T_{\infty}(B_w \otimes \omega_w) T_{\infty}^{-1}Q_2^{-1} +  R_0^{-1}(A_z^0 \otimes \omega_z )R_0 ) & \text{on } U_{\infty},
\end{cases}
\end{equation*}
which has bundle type $\cO(1) \oplus \cO(-2)$.

Secondly, we determine the remainder parameters of $B_w$ so that the residue matrix of $\nabla_0+ \Phi_{X}$ at $t_i$ has 
the eigenvalues $(\nu_i,-\nu_i)$ (resp. $(\nu_i,1-\nu_i)$) for $i=1,\ldots,4$ (resp. $i=5$).
We take eigenvectors of the residue matrices of $\Phi_{X}$ at $t_i$ as follows:
\begin{equation*}
\begin{cases}
(v_{\nu_{i}'^+}^{1},v_{\nu_{i}'^+}^{2}) = (f_i^+(q_1,p_1,\lambda), (q_1-q_2)(q_1-t_i)  (q_2-t_i) ) & \text{ associated to $\nu_i'$} \\
(v_{\nu_{i}'^-}^{1},v_{\nu_{i}'^-}^{2}) = (f_i^-(q_1,p_1,\lambda), (q_1-q_2)(q_1-t_i)  (q_2-t_i) ) & \text{ associated to $-\nu_i'$}.
\end{cases}
\end{equation*}
Here, we put
$f_i^{\epsilon}(q_1,p_1,\lambda) =p_1(q_2 - t_i)  +p_2 (q_1 - t_i)  - \epsilon \nu_i'  (q_1 - q_2) (t_i - t_j) (t_i - t_k) (t_i - t_l)$
where $\epsilon \in \{+,-\}$.
Set
\begin{equation*}
T_{t_i}^{\epsilon}=\left(
\begin{array}{ll}
1 & f_i^{\epsilon}(q_1,p_1,\lambda)  \\
0 & (q_1-q_2)(q_1-t_i)  (q_2-t_i)  
\end{array}
\right)  \quad i=1,\ldots,4, \ \epsilon \in \{+,-\}.
\end{equation*}
Then we have
\begin{equation*}
(T_{t_i}^{\epsilon})^{-1} (\res_{t_i}(Q_1^{-1}P^{-1}(A_z^0\otimes \omega_z) P Q_1))T_{t_i}^{\epsilon} = 
\left(
\begin{array}{ll}
-\epsilon\nu'_i & *  \\
0 & \epsilon\nu'_i  
\end{array}
\right).
\end{equation*}

We fix a tuple $\boldsymbol{\epsilon}=(\epsilon_1,\ldots,\epsilon_4)$ where $\epsilon_i \in \{ +,-\}$.
For the tuple $\boldsymbol{\epsilon}$,
we solve the following four linear equations
\begin{equation*}
( (T_{t_i}^{\epsilon_i})^{-1}(\res_{t_i} (R_{0,w}^{-1} d R_{0,w} + R_{0,w}^{-1} T_{\infty}(B_w \otimes \omega_w) T_{\infty}^{-1} R_{0,w})) T_{t_i}^{\epsilon_i})_{21} =0\quad i=1,\ldots,4
\end{equation*}
where $A_{ij}$ is the $(i,j)$-entry of a $(2\times 2)$-matrix $A$.
Explicitly, we can describe these equations as follows:
\begin{equation}\label{flag equation}
\begin{aligned}
&(q_1-q_2)(q_1 - t_i)^2 (q_2 - t_i)^2 \bigl( t_i^3 d_1+  t_i^2 d_2+ t_i d_3 +  d_4 \bigr) \\
&\quad +2 (q_1 - t_i)  (q_2 - t_i)  g_i(q_1,p_1,\lambda ) \bigl( t_i^3 f_0+ t_i^2  f_1 + t_i f_2 +f_3 \bigr)\\
&\qquad - ( g_i^{\epsilon_i}(q_1,p_1,\lambda ))^2 \bigl( t_i^2 e_0+t_i e_1 +e_2 \bigr) =0  \quad i=1,\ldots,4
\end{aligned}
\end{equation}
where we put
\begin{equation*}
\begin{aligned}
g_i^{\epsilon_i}(q_1,p_1,\lambda ) :=&\ p_1(q_2 - t_i)  + p_2(q_1 - t_i) -\nu_5' t_i  (q_1 - q_2) (q_1 - t_i) (q_2 - t_i) +\epsilon_i \nu_{i}' (q_1 - q_2) (t_i - t_j) (t_i - t_k)(t_i- t_l)\\
   =&\ p_1 (q_1 -t_i) + p_1 (q_2 -  t_i) - (q_1 - q_2) ( (q_1-t_i)\lambda + \nu_5' t_i (q_1-t_i) (q_2 - t_i) - \epsilon_i \nu_i' (t_i-t_j) (t_i - t_k) (t_i-t_l) ).
\end{aligned}
\end{equation*}
Then the solution of the equations (\ref{flag equation}) is the following
\begin{equation}\label{d1d2d3d4}
\begin{pmatrix}
d_1\\
d_2\\
d_3\\
d_4
\end{pmatrix}=
\begin{pmatrix}
0 & 0 & 0 & 1\\
 1  & 1 & 1 & 1  \\
x_1^3 & x_1^2 & x_1 & 1\\
x_2^3 & x_2^2 & x_2 & 1
\end{pmatrix}^{-1}
\begin{pmatrix}
\frac{h_1^{\epsilon_1}(q_1,p_1,\lambda)}{q_1^2 q_2^2}\\
\frac{h_2^{\epsilon_2}(q_1,p_1,\lambda)}{(q_1-1)^2 (q_2-1)^2}\\
\frac{h_3^{\epsilon_3}(q_1,p_1,\lambda)}{(q_1-x_1)^2 (q_2-x_1)^2}\\
\frac{h_4^{\epsilon_4}(q_1,p_1,\lambda)}{(q_1-x_2)^2 (q_2-x_2)^2}
\end{pmatrix}
\end{equation}
where
\begin{equation*}
\begin{aligned}
h_i^{\epsilon_i}(q_1,p_1,\lambda)=& \frac{1}{2} (q_1 - t_i)  \bigl( -2 e_1 p_1 -  2 e_0 p_1 (q_1 + t_i) + (q_1 -t_j) (q_1 -t_k) (q_1 - t_l) \bigr)\Bigl(p_1 (q_1 - t_i) + p_1 (q_2 - t_i) \\
&- 2  (q_1 - t_i) \bigl( (q_1 - t_i)\lambda +  \nu_5' t_i (q_1 - t_i) (q_2 - t_i)  - \epsilon_i \nu_i'  (t_i-t_j) (t_i-t_k) (t_i- t_l) \bigr) \\
&+ \frac{q_1 - q_2}{p_1} \bigl( (q_1 - t_i)\lambda +  \nu_5' t_i (q_1 - t_i) (q_2 - t_i)  - \epsilon_i\nu_i' (t_i-t_j) (t_i-t_k) (t_i- t_l) \bigr)^2 \Bigr).
\end{aligned}
\end{equation*}
We consider the case $q_1 = t_i$ for some $i$.
For fixed $\boldsymbol{\epsilon}=(\epsilon_1,\ldots,\epsilon_5)$, $\epsilon_i \in \{ +,- \}$, 
the domain of the function $\frac{h_i^{\epsilon_i}(q_1,p_1,\lambda)}{(q_1-t_i)^2 (q_2-t_i)^2}$ is extended to $q_1=t_i$, $p_1=\epsilon_i \hat{\nu}_i$ and $\lambda=v^{i,\epsilon_i}_1$.
Here, $v^{i,\epsilon}_1$ is the blowing-up parameter of $\cK_5' \ra \cK_5$ at $(t_i,\epsilon_i \hat{\nu}_i )$.
Therefore we can extend the family $\nabla_0$ to 
$q_1=t_i$, $p_1=\epsilon_i \hat{\nu}_i$ and $\lambda=v^{i,\epsilon_i}_1$ when we substitute the solution $d_1,\ldots,d_4$ associated to $\boldsymbol{\epsilon}$.

We compute the eigenvalue of the residues of $\nabla_0$ at $t_i$ for $i=1,\ldots,5$.
We put
\begin{equation}\label{eigen alpha}
\begin{aligned}
\alpha_i^{\epsilon_i}(q_1,q_2,p_1,\lambda) 
:=&\ ( {T_{t_i}^{\epsilon}}^{-1}(\res_{t_i} (R_{0,w}^{-1} d R_{0,w} + R_{0,w}^{-1} T_{\infty}(B_w \otimes \omega_w) T_{\infty}^{-1} R_{0,w}))) T_{t_i}^{\epsilon})_{22}\\
= &\ \frac{\beta_i^{(1)}(q_1,p_1,\lambda) \beta^{(2)}_i(\epsilon_i,q_1,p_1,\lambda)}{2 p_1 (q_2-t_i) (t_i -t_j) (t_i- t_k) (t_i-t_l)} \quad i=1,\ldots,4,
\end{aligned}
\end{equation}
where
\begin{equation*}
\begin{aligned}
\beta^{(1)}_i(q_1,p_1,\lambda) &:= -2 e_1 p_1 - 2 e_0 p_1 (q_1 + t_i) + (q_1-t_j) (q_1-t_k)  (q_1 - t_l), \\
\beta^{(2)}_i(\epsilon_i, q_1,p_1,\lambda) &:= (p_1 (q_1 - t_i) - (q_1 - q_2) ((q_1 - t_i) \lambda +\nu_5' t_i (q_1 - t_i) (q_2 - t_i) - \epsilon_i \nu_i' (t_i-t_j) (t_i - t_k) (t_i-t_l)  )).
\end{aligned}
\end{equation*}
Then we have
\begin{equation*}
(T_{t_i}^{\epsilon_i})^{-1} (\res_{t_i}((\nabla_0 )_{U_0})T_{t_i}^{\epsilon_i} = 
\left(
\begin{array}{ll}
-\alpha_i^{\epsilon_i}(q_1,q_2,p_1,\lambda) & *  \\
0 & \alpha_i^{\epsilon_i}(q_1,q_2,p_1,\lambda) 
\end{array}
\right)
\end{equation*}
for $i=1,\ldots,4$.
Here we put $(\nabla_0)_{U_0}= R_{0,w}^{-1} d R_{0,w} + R_{0,w}^{-1} T_{\infty}(B_w \otimes \omega_w) T_{\infty}^{-1} R_{0,w}$.

We define the family of connections $\nabla_X$ as $\nabla_0+\Phi_{X}$.
The limit $\lim_{q_2\ra q_1} \nabla_X$ has bundle type $\cO(1)\oplus \cO(-2)$.
The eigenvalues of the residues of $\nabla_X$ at $0,1,x_1,x_2$, and $\infty$ are the following
\begin{center}
  \begin{tabular}{|c|c|c|c|c|} \hline
    $\res_{0} \nabla_X$ & $\res_{1} \nabla_X$ 
& $\res_{t_1} \nabla_X$ & $\res_{t_2} \nabla_X$ & $\res_{\infty} \nabla_X$ \\ \hline 
   $-\alpha_1^{\epsilon_1} +\epsilon_1 \nu_1'$  
&  $-\alpha_2^{\epsilon_2} + \epsilon_2 \nu_2'$ 
&$-\alpha_3^{\epsilon_3} + \epsilon_3 \nu_3'$
&$-\alpha_4^{\epsilon_4} + \epsilon_4 \nu_4'$
&$-\frac{1}{2} + \nu_5'$  \\ \hline
   $ \alpha_1 ^{\epsilon_1}-\epsilon_1 \nu_1' $  
&  $ \alpha_2 ^{\epsilon_2}- \epsilon_2 \nu_2' $ 
&$ \alpha_3 ^{\epsilon_3}- \epsilon_3 \nu_3' $
&$\alpha_4 ^{\epsilon_4}- \epsilon_4 \nu_4' $
&$ \frac{3}{2}- \nu_5' $  \\ \hline
  \end{tabular}
\end{center}
where we put $\alpha^{\boldsymbol{\epsilon}}_i:=\alpha^{\boldsymbol{\epsilon}}_i(q_1,q_2,p_1,\lambda)$.
For the fixed tuple $(\nu_1,\ldots,\nu_5)$, we determine the eigenvalues $(\nu'_1,\ldots,\nu'_5)$ of Higgs field $\Phi_{X}$ as follows.
We solve the following equations for $\nu_i'$:
\begin{equation}\label{condition of eigen}
-\frac{1}{2} + \nu_5'= \nu_5, \text{ and } -\alpha_i^{\epsilon_i} + \epsilon_i \nu_i'=\nu_i \text{ for $i=1,\ldots,4$}.
\end{equation}
Then $\nu_i'$ is determined by $\nu_i$ and $\nu_5$ (when at least $q_2$ is close to $q_1$ enough).
Then the eigenvalues of residues of the family $\nabla_X$ at $t_i$ ($i=1,\ldots,5$) are $\nu_i^{\pm}$.
Let $\widetilde{A}_{U_0}$ be the connection matrix
\begin{equation*}
R_{0,w}^{-1} d R_{0,w} + R_{0,w}^{-1} T_{\infty}(B_w \otimes \omega_w) T_{\infty}^{-1} R_{0,w}+
 Q_1^{-1}P^{-1}(A_z^0\otimes \omega_z) P Q_1 
\end{equation*}
(see (\ref{Def of jumping family})).
Then, we obtain the family (\ref{jump family descrip conn}), which satisfies the assertion of Proposition \ref{Jump descrip Prop} by the construction.
\end{proof}

\begin{Rem}
The connection $\nabla_X$ is parametrized by $e_0$ and $e_1$.
We can take $e_0,e_1 \in \mathbb{C}$ freely.
\end{Rem}

\subsection{The apparent singularities and the dual parameters of $\nabla_X$}\label{lim of app and dual}
Let $C_{\infty}$ be the $\infty$-section of the Hirzebruch surface of $n-2$.
Let $\phi \colon \widehat{M} \dashrightarrow \Hilb^2(\widetilde{\cK}_n')$ be the birational map constructed by the apparent singularities and the dual parameters.
For $n=5$, by analyzing the behavior of the apparent singularities and their duals of $\nabla_X$ as $q_2 \ra q_1$, we obtain the following 
\begin{Thm}\label{main theorem conn}
\textit{
By taking a sequence of blowing-ups $\widetilde{\Hilb}^2(\widetilde{\cK}_5'\cup C_{\infty}) \ra \Hilb^2(\widetilde{\cK}_5'\cup C_{\infty})$,
we have the injective map $\tilde{\phi} \colon \widehat{M} \rightarrow \widetilde{\Hilb}^2(\widetilde{\cK}_5'\cup C_{\infty})$.
The moduli space $\widehat{M}$ is biregular to its image $\mathrm{Image}\, \tilde{\phi}(\widehat{M}) \subset \widetilde{\Hilb}^2(\widetilde{\cK}_5'\cup C_{\infty})$.}
\end{Thm}

\begin{proof}
First, we construct a map $\widehat{M} \ra \Hilb^2(\widetilde{\cK}_5'\cup C_{\infty})$ by the birational map $\phi$.
The apparent singularities of $\nabla_X$ is the zero of the polynomial
\begin{equation*}
\begin{aligned}
F_{21}(z)=&\ (1 + e_0 (q_1 - q_2))z^2 - (q_1 +q_2 + e_1( q_2 - q_1) ) z\\
&\qquad \qquad  + q_1 q_2  + \frac{(q_1 - q_2)q_1}{2 p_1} (-2 e_1 p_1 - 2 e_0 p_1 q_1 + (q_1-1) (q_1 - x_1) (q_1 - x_2)).
\end{aligned}
\end{equation*}
Let $q_1'$ and $q_2'$ be the solution of the equation $F_{21}(z)=0$:
\begin{equation*}
q_1'=q_1+a_1(q_2-q_1)- \frac{\sqrt{a_3(q_2-q_1)}}{a_2}, \quad q_2'=q_1+a_1(q_2-q_1)+ \frac{\sqrt{a_3(q_2-q_1)}}{a_2}
\end{equation*}
where we set
\begin{equation*}
\begin{aligned}
a_1&:=\frac{1 + e_1 + 2 e_0 q_1}{2 - 2 e_0 (q_2-q_1)}, \quad a_2:=2 p_1 (-1 + e_0 (q_2-q_1)), \\
a_3&:=p_1 (-p_1 (1 + e_1 + 2 e_0 q_1)^2 (q_1 - q_2) + 2 q_1(q_1-1)(1 + e_0 (q_1 - q_2)) (q_1 - x_1) (q_1 - x_2)).
\end{aligned}
\end{equation*}
Note that $\lim_{q_2\ra q_1} q_1'=\lim_{q_2\ra q_1} q_2'=q_1$.
The dual parameter of $q_i'$ is $p_i'=F_{11}(q_i')$ for $i=1,2$.
If we consider $\lim_{q_2\ra q_1}F_{11}(q_i')$, then the limit diverges
(Figure \ref{fig3}).
\begin{figure}[htbp]
 \begin{center}
  \includegraphics[width=55mm]{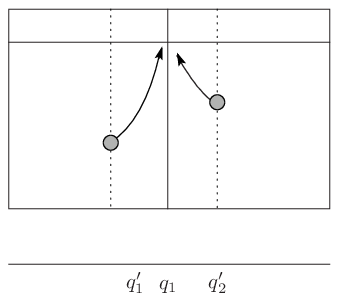}
 \end{center}
 \caption{The apparent singularities and the dual parameters}
 \label{fig3}
\end{figure}
Let $s$ be the parameter such that $\bar{p}'_2-\bar{p}'_1 =s(q_2'-q_1')$ where $\bar{p}'_i=1/p_i$ for $i=1,2$.
By explicit computation, we have
\begin{equation*}
\lim_{q_2\ra q_1} s=\frac{1}{ q_1 (q_1-1) (q_1-x_1) (q_1-x_2)}.
\end{equation*}
Then we can extend the map $X \ra \mathrm{Hilb}^2(\widetilde{\cK}_5')$ to $\widehat{X} \ra \mathrm{Hilb}^2(\widetilde{\cK}_5'\cup C_{\infty})$ by 
\begin{equation}
\widehat{X} \setminus X \ni ((q_1,p_1),(q_1,p_1), \lambda) \longmapsto ( \{ (q_1 ,\infty ), (q_1 ,\infty ) \}, \frac{1}{ q_1 (q_1-1) (q_1-x_1) (q_1-x_2)} ) \in \mathrm{Hilb}^2(\widetilde{\cK}_5'\cup C_{\infty}).
\end{equation}
We can define a map $\phi \colon \widehat{M} \ra \Hilb^2(\widetilde{\cK}_5'\cup C_{\infty})$
so that $\phi|_{M^0}$ is the map in Theorem \ref{from M0con to Hilb} and the diagram
\begin{equation*}
\xymatrix@=15pt{
  \widehat{X} \ar[d] \ar[dr]   & \\
\widehat{M} \ar[r]^-{\phi}  & \Hilb^2(\widetilde{\cK}_5'\cup C_{\infty}).
}
\end{equation*}
is commutative.

Now, we construct a sequence of blowing-ups $\widetilde{\Hilb}^2(\widetilde{\cK}_5'\cup C_{\infty})$ of $\Hilb^2(\widetilde{\cK}_5'\cup C_{\infty})$ such that
\begin{equation*}
\xymatrix@=15pt{
 & \widetilde{\Hilb}^2(\widetilde{\cK}_5' \cup C_{\infty})  \ar[d] \\
\widehat{M} \ar[r]_-{\phi} \ar[ur]^-{\tilde{\phi}} & \Hilb^2(\widetilde{\cK}_5'\cup C_{\infty}) 
}
\end{equation*}
where $\tilde{\phi}$ is injective.
Set
$\bar{p}'_i:=\frac{1}{p_i'} $.
Let $((q_1',\bar{p}_1'), (q_2',\bar{p}_2'))$ be coordinates on 
\begin{equation*}
(\widetilde{\cK}_5'\cup C_{\infty})\times (\widetilde{\cK}_5' \cup C_{\infty})\setminus (C_0\cup\pi^{-1}(\infty)) \times (C_0\cup\pi^{-1}(\infty)).
\end{equation*}
First, we take the blowing-up of $(\widetilde{\cK}_5'\cup C_{\infty})\times (\widetilde{\cK}_5' \cup C_{\infty})$ along the ideal $\cI_1= (q_2'-q_1', \bar{p}'_2-\bar{p}'_1)$, denoted by $\mathrm{Bl}_1$.
Note that $\Hilb^2(\widetilde{\cK}_5'\cup C_{\infty}) = \mathrm{Bl}_1/\mathfrak{S}_2$.
We defined the blowing-up parameter $s$ as $\bar{p}'_2-\bar{p}'_1 =s(q_2'-q_1')$, and
we obtained
\begin{equation*}
\lim_{q_2\ra q_1} s=\frac{1}{ q_1 (q_1-1) (q_1-x_1) (q_1-x_2)}.
\end{equation*}
Second, we take the blowing-up of $\mathrm{Bl}_1$ along the ideal
\begin{equation*}
\cI_2= \left( q_2'-q_1', \bar{p}'_2-\bar{p}'_1, s-\frac{1}{ q (q-1) (q-x_1) (q-x_2)} , \bar{p}'_2+\bar{p}'_1 \right) \quad \text{where $q:=\frac{q_1'+q_2'}{2}$},
\end{equation*}
denoted by $\mathrm{Bl}_2$.
We define the blowing-up parameters $t_1$, $t_2$ as 
\begin{equation*}
 s-\frac{1}{ q (q-1) (q-x_1) (q-x_2)}= t_1 (q_2'-q_1'), \quad  \bar{p}'_2+\bar{p}'_1 =t_2(q_2'-q_1').
\end{equation*}
By explicit calculations, we have $\lim_{q_2\ra q_1} t_1=\lim_{q_2\ra q_1} t_2=0$.
Third, we take the blowing-up of $\mathrm{Bl}_2$ along the ideal 
\begin{equation*}
\cI_3= \left(q_2'-q_1', \bar{p}'_2-\bar{p}'_1, s-\frac{1}{ q (q-1) (q-x_1) (q-x_2)} , \bar{p}'_2+\bar{p}'_1,  t_1 ,t_2 \right),
\end{equation*}
denoted by $\mathrm{Bl}_3$.
We define the blowing-up parameters $u_1$, $u_2$ as 
$ t_1= u_1 (q_2'-q_1'), \quad t_2= u_2 (q_2'-q_1')$.
By explicit calculations, we have
\begin{equation}\label{lim u}
\begin{aligned}
\lim_{q_2\ra q_1} u_1=
&\frac{-\lambda  }{4 q_1^2 (q_1-1)^2 (q_1-x_1)^2 (q_1-x_2)^2}  +\frac{ u_{1}^{(1)}(q_1) p_1}{4 q_1^3 (q_1-1)^3 (q_1-x_1)^3 (q_1-x_2)^3} \\
& \quad + \frac{u_{1}^{(0)}(q_1)}{16 q_1^3 (q_1-1)^3 (q_1-x_1)^3 (q_1-x_2)^3} ,\\
\lim_{q_2\ra q_1} u_2=&\frac{-4 q_1^3+ 3 q_1^2 (1 + x_1 + x_2) - 2 q_1 (x_1 + x_2 + x_1 x_2)+ x_1 x_2}{ 4 q_1^2 (q_1-1)^2 (q_1-x_1)^2 (q_1-x_2)^2}
\end{aligned}
\end{equation}
where $u_{1}^{(1)}(q_1)$ and $u_{1}^{(0)}(q_1)$ are polynomials in $q_1$, which are independent of $\lambda$ and $p_1$.
Fourth, we take the blowing-up of $\mathrm{Bl}_3$ along the ideal
\begin{equation*}
\begin{aligned}
\cI_4= \Big(q_2'-q_1', &\bar{p}'_2-\bar{p}'_1, s-\frac{1}{ q (q-1) (q-x_1) (q-x_2)} , \bar{p}'_2+\bar{p}'_1,  t_1 ,t_2,\\
&u_2 -\frac{-4 q^3+ 3 q^2 (1 + x_1 + x_2) - 2 q (x_1 + x_2 + x_1 x_2)+ x_1 x_2}{ 4 q^2 (q-1)^2 (q-x_1)^2 (q-x_2)^2} \Big),
\end{aligned}
\end{equation*}
denoted by $\mathrm{Bl}_4$.
We put blowing-up parameters $v$ as follows:
\begin{equation*}
 u_2-\frac{-4 q^3+ 3 q^2 (1 + x_1 + x_2) - 2 q (x_1 + x_2 + x_1 x_2)+ x_1 x_2}{ 4 q^2 (q-1)^2 (q-x_1)^2 (q-x_2)^2}= v (q_2'-q_1').
\end{equation*}
By explicit calculations, we have $\lim_{q_2\ra q_1}v=0$.
Finally, we take the blowing-up of $\mathrm{Bl}_4$ along the ideal 
\begin{equation*}
\begin{aligned}
\cI_5= \Big(q_2'-q_1', &\bar{p}'_2-\bar{p}'_1, s-\frac{1}{ q (q-1) (q-x_1) (q-x_2)} , \bar{p}'_2+\bar{p}'_1,  t_1 ,t_2,\\
&u_2 -\frac{-4 q^3+ 3 q^2 (1 + x_1 + x_2) - 2 q (x_1 + x_2 + x_1 x_2)+ x_1 x_2}{ 4 q^2 (q-1)^2 (q-x_1)^2 (q-x_2)^2}, v \Big),
\end{aligned}
\end{equation*}
 denoted by $\mathrm{Bl}_5$.
We define the blowing-up parameter $w$ as $v= w (q_2'-q_1')$.
By explicit calculations, we have
\begin{equation}\label{lim w}
\begin{aligned}
\lim_{q_2\ra q_1} w=
&\frac{\lambda}{ 8 q_1^3 (q_1-1)^3 (q_1-x_1)^3 (q_1-x_2)^3}(4 q_1^3- 3 q_1^2 (1 + x_1 + x_2) + 2 q_1 (x_1 + x_2 + x_1 x_2)- x_1 x_2)\\
&+\frac{-w^{(1)}(q_1) p_1}{8 q_1^4 (q_1-1)^4 (q_1 - x_1)^4 (q_1 - x_2)^4}  +\frac{w^{(0)}(q_1)}{64 q_1^4 (q_1-1)^4 (q_1 - x_1)^4 (q_1 - x_2)^4}
\end{aligned}
\end{equation}
where $w_{1}^{(1)}(q_1)$ and $w_{1}^{(0)}(q_1)$ are polynomials in $q_1$, which are independent of $\lambda$ and $p_1$.
We define the blowing-ups $\widetilde{\Hilb}^2(\widetilde{\cK}_5')$ of $\Hilb^2(\widetilde{\cK}_5')$ as $\widetilde{\Hilb}^2(\widetilde{\cK}_5') = \mathrm{Bl}_5 / \mathfrak{S}_2$.
Here, $\mathfrak{S}_2$ is the symmetric group, which acts on $\mathrm{Bl}_5$ naturally.

We consider the family $\lim_{q_2\ra q_1} \nabla_X$, which has bundle type $\cO(1)\oplus \cO(-2)$.
The family is parameterized by $(\lambda, p_1,q_1) \in \mathbb{C}^3$.
We fix the parameter $q_1$.
We can regard the parameters $(\lambda,p_1)$ as coordinates of $M^1$, which is isomorphic to $\mathbb{C}^2$.
The family $\lim_{q_2\ra q_1} \nabla_X$ where $q_1$ is fixed is a universal family of $M^1$.
The theorem follows form (\ref{lim u}) and (\ref{lim w}).
\end{proof}


\begin{thebibliography}{99}

\bibitem{AL}
D.\ Arinkin, S.\ Lysenko, 
{\it On the moduli of $SL(2)$-bundles with connections on $\bP^1\setminus \{ x_1,\ldots,x_4 \}$}, 
Internat. Math. Res. Notices (1997), no. {\bf 19}, 983--999.

\bibitem{BNR}
A.\ Beauville, M.\ S.\ Narasimhan, S.\ Ramanan, 
{\it Spectral curves and the generalised theta divisor}.
J. Reine Angew. Math. {\bf 398} (1989), 169--179.

\bibitem{Drinfeld}
V. G. Drinfeld, 
{\it Elliptic modules and their applications to the Langlands and to
the Peterson conjectures for $\mathrm{GL}(2)$ over functional field} (in Russian), 
Ph. D. Thesis, Moscow State University, 1977.


\bibitem{DM}
B.\ Dubrovin, M.\ Mazzocco,
{\it Canonical structure and symmetries of the Schlesinger equations}. 
Comm. Math. Phys. {\bf 271} (2007), no. 2, 289--373. 

\bibitem{Inaba}
M.-A.\ Inaba,
{\it Moduli of parabolic connections on curves and the Riemann-Hilbert correspondence}. 
J. Algebraic Geom. {\bf 22} (2013), no. 3, 407--480.

\bibitem{IIS}
M.\ Inaba, K.\ Iwasaki, M.-H.\ Saito, 
{\it Moduli of stable parabolic connections, Riemann-
Hilbert correspondence and geometry of Painlev\'e equation of type VI. I },
Publ. Res. Inst. Math. Sci. (2006), no. {\bf 4}, 987-1089.


\bibitem{IIS2}
M.\ Inaba, K.\ Iwasaki, M.-H.\ Saito, 
{\it Moduli of stable parabolic connections, Riemann-Hilbert correspondence and geometry of Painlev\'e equation of type VI. II}.
Moduli spaces and arithmetic geometry, 387-432, Adv. Stud. Pure Math., {\bf 45}, Math. Soc. Japan, Tokyo, 2006.


\bibitem{IS}
M.\ Inaba, M.-H.\ Saito, 
{\it Moduli of regular singular parabolic connections of spectral type on smooth projective curves}.
In preparation.


\bibitem{LS}
F.\ Loray, M.-H.\ Saito,
{\it Lagrangian fibrations in duality on moduli spaces of rank 2 logarithmic connections over the projective line.} 
Internat. Math. Res. Notices (2015), no. {\bf 4}, 995--1043.



\bibitem{Naka}
H.\ Nakajima,
{\it Lectures on Hilbert schemes of points on surfaces}. University Lecture Series, {\bf18}. American Mathematical Society, Providence, RI, 1999.


\bibitem{Okamoto}
K. Okamoto,
{\it Isomonodromic deformation and Painlev\'e equations, and the Garnier system.}
J. Fac. Sci. Univ. Tokyo Sect. IA Math. {\bf 33} (1986), no. 3, 575--618.

\bibitem{Obl}
S.\ Oblezin, 
{\it Isomonodromic deformations of $\mathfrak{sl}(2)$ Fuchsian systems on the Riemann sphere.}
Mosc. Math. J. 5 (2005), no. {\bf 2}, 415--441, 494--495.


\bibitem{SS}
M.-H.\ Saito, S.\ Szabo,
{\it Apparent singularities and canonical coordinates of moduli spaces of parabolic connections and parabolic Higgs bundles}.
in preparation.


\bibitem{Sim}
C.\ Simpson,
{\it Harmonic bundles on noncompact curves}. 
J. Amer. Math. Soc. 3 (1990), no. {\bf 3}, 713--770

\bibitem{Skl}
E.\ K.\ Sklyanin, 
{\it Separation of variables---new trends}, 
Progr. Theoret. Phys. Suppl. (1995),
no. 118, 35--60.

\bibitem{Yoshi}
H.\ Yoshida,
{\it On the structure of strata of the moduli space of parabolic connections with 5-regular singular points on $\bP^1$},
(Japanese), master thesis, Kobe university, (2015).


\end{thebibliography}
\end{document}